\documentclass[12pt]{iopart}

  \expandafter\let\csname equation*\endcsname\relax
  \expandafter\let\csname endequation*\endcsname\relax
\usepackage{amsmath}
\usepackage{amsfonts}
\usepackage{iopams}
\usepackage{cite}
\usepackage{rotating}
\usepackage{graphicx}
\usepackage{chngcntr}
\counterwithin{figure}{section}

\newtheorem{theorem}{Theorem}

\newtheorem{remark}{Remark}
\newtheorem{proposition}{Proposition}
\newtheorem{observation}{Observation}

\usepackage{color}

\begin{document}


\title[A Family of Periodic Orbits in the Three-Dimensional Lunar Problem]{A Family of Periodic Orbits in the Three-Dimensional Lunar Problem}

\author{Edward Belbruno$^1$, Urs Frauenfelder$^2$, Otto van Koert$^3$}

\address{$^1$ Department of Astrophysical Sciences, Princeton University, Princeton, NJ, USA; Department of Mathematics,  Yeshiva University,  New York, NY, USA }
\address{$^2$ Institut fur Mathematik, Augsburg University, Augsburg, Germany}
\address{$^3$ Department of Mathematical sciences, Seoul National University, Seoul, South Korea}

\eads{\mailto{belbruno@princeton.edu}, \mailto{urs.frauenfelder@math.uni-augsburg.de},
\mailto{okoert@snu.ac.kr}
}

\begin{abstract}
A family of periodic orbits is proven to exist in the spatial lunar problem that are continuations of a family of consecutive collision orbits, perpendicular to the primary orbit plane. This family emanates from all but two energy values. The orbits are numerically explored. 
The global properties and geometry of the family is studied. 
\end{abstract}



\section{Introduction} 
\label{sec:intro} 

We consider the three-dimensional circular restricted three-body problem. This models the three-dimensional motion of a particle, $P_0$, of zero mass in the Newtonian gravitational field generated by two particles, $P_1, P_2$ of respective positive masses, $m_1, m_2$, in a mutual uniform circular motion. It is assumed that $m_1$ is much larger than $m_2$. This problem is studied in a rotating coordinate system that rotates with the same constant frequency, $\omega$ of the circular motion of $P_1, P_2$, so that in this system $P_1$ and $P_2$ are fixed. 
Because $m_1$ is much larger than $m_2$, we refer to $P_1$ as the Earth and $P_2$ as the Moon, for convenience. 

When $P_0$ moves about the larger particle, $P_1$, the motion of $P_0$ can be completely understood  if, for example, $P_0$ is restricted to the two-dimensional plane of motion of $P_1, P_2$. In this case, with $m_2 = 0$, assume that $P_0$ has precessing elliptic motion, of elliptic frequency $\omega^*$ about $P_1$, precessing with frequency $\omega$. Then the Kolmogorov-Arnold-Moser(KAM) Theorem proves that this precessing motion persists if $m_2$ is sufficiently small and if $\omega$ and $\omega^*$ are sufficiently noncommensurate. Otherwise, the motion is chaotic due to heteroclinic dynamics. That is, invariant KAM tori foliate the phase space. The motion of $P_2$ is proven to be stable \cite{SiegelMoser:1971}. When the initial elliptic motion of $P_0$ is not in the same plane as the $P_0, P_1$ then under similar assumptions although KAM tori can be proven to exist, but stability cannot be guaranteed. 

In this paper we study the three-dimensional motion of $P_0$ about $P_2$. This is referred to as the three-dimensional, or spatial, lunar problem. Relatively little is proven in general about the motion of $P_0$ unless the initial motion starts infinitely close to $P_2$.   The proof of existence of KAM tori in the three-dimensional lunar problem was obtained by M. Kummer under the assumption that the initial motion of $P_0$ lies infinitely near to  
$P_2$ \cite{Kummer:1983}. \footnote{Kummer proved the existence of KAM tori in the planar lunar problem sufficiently near to $P_2$ \cite{Kummer:1971}. This also proves the stability of the Hill periodic orbits.
 \cite{Hill:1878}.}  

The main result of this paper is to prove the existence of a special family of periodic orbits about $P_2$, nearly perpendicular to the primary orbit plane.  More precisely, if we normalize $m_1 = 1 - \mu , m_2 = \mu$, then in the case of $\mu = 0$ there exists a family of periodic orbits on the $z-$axis through $P_2$, so perpendicular to the $P_1, P_2$ plane, parameterized by their energy $h$. This family consists of consecutive collision orbits: Starting at collision at $P_2$, they extend up the $z-$axis to a maximal distance $d=d(h)$, then fall back to $P_2$, and periodically repeat this oscillation, where $d$ can have any positive value. We label these as $\phi^*(t,h)$. These orbits have period $T^*(h)$. We prove that $\phi^*(t,h)$ smoothly continues as a function of $\mu$, sufficiently small, into a unique periodic orbit $\phi(t,\mu)$ of period $T = T^* + \mathcal{O}(\mu)$, on the associated Jacobi energy surface, provided a non-degeneracy condition holds. This condition is satisfied for every energy value except two. The resulting family periodic orbits is labeled, $ \mathcal{F}(h, \mu)$.

The method of proof of $ \mathcal{F}(h, \mu)$ is to make use of the proof of existence of an analogous family of orbits about the primary mass $P_1$ \cite{Belbruno:1981a}. A three-dimensional regularization defined first in \cite{Belbruno:1981a} is performed. This uses a fractional linear M\"obius transformation which can be represented elegantly using a Jordan algebra. We also exploit symmetry properties of the lunar problem.

The resulting family of orbits is numerically investigated and has interesting properties. The properties are analogous to those in \cite{Belbruno:1981b} for very negative energy, but differ markedly for larger energy.

Geometric and global properties of the larger family of periodic orbits connecting those found here and those in \cite{Belbruno:1981b} are described using ideas from symplectic geometry \cite{Frauenfelder:2018}. In particular, we study the extension of  $\mathcal{F}(h, \mu)$  as a function of $\mu$ by using a construction due to P. Rabinowitz \cite{Rabinowitz}.

The main theorem for this paper is stated as Theorem 1 in Section \ref{sec:Result}. The proof is done in 
Section~\ref{sec:Proof}. In Section \ref{sec:Numerical}, we describe numerical results. Global and geometric properties are described in Section~\ref{sec:Geometry} and summarized in Theorems 2,3.

\section{Spatial Lunar Problem and Main Result}
\label{sec:Result} 

The three-dimensional restricted three-body problem in a rotating coordinate system with coordinates $q = (q_1, q_2, q_3) \in \mathbb{R}^3$ and momenta, $p = (p_1, p_2, p_3) \in \mathbb{R}^3$, for the motion of the zero mass particle $P_0$ is defined by the Hamiltonian system,
\begin{equation}
H =  \frac{1}{2}|p|^2-\frac{\mu}{|q-m|}-\frac{1-\mu}{|q-e|}+ \omega(q_1 p_2-q_2p_1) ,
\label{eq:H}
\end{equation}
\begin{equation}
\dot{q} = H_p , \quad \dot{p} = - H_q,
\label{eq:H-DE}
\end{equation}
$^. \equiv d/dt $, $t \in  \mathbb{R}^1$ is time, $H_p \equiv \partial H/\partial p$,
where the masses of $P_1, P_2$ are normalized to be $m_1 = 1 - \mu, \quad m_2 = \mu$, respectively, $\mu \in (0,1]$. The center of mass is placed at the origin. $P_1, P_2$ are fixed on the $q_1$-axis at the respective locations,
$ e=(-\mu,0,0), \quad m=(1-\mu,0,0)$, where $e, m$ denotes Earth, Moon, respectively. $\omega$ is the frequency of the rotating frame that is normalized to be $1$, where we use $\omega$ for generality. A solution $\phi = \phi(t) \in \mathbb{R}^6$ of \eqref{eq:H-DE} lies on the $5$-dimensional energy surface
\begin{equation}
\Sigma_H = \{ (q,p) \in \mathbb{R}^6 | H(q,p) = -h \}, 
\label{eq:H-Surface}
\end{equation}
where $ h  \in \mathbb{R}^1$ is a constant. (In nonsymplectic coordinates, $H$ can be written in another form called the Jacobi integral.)
\medskip

\noindent
{\bf Definition.}  The {\em spatial lunar problem} is defined by viewing the motion of $P_0$ from \eqref{eq:H-DE} to lie near $P_2$ and assuming that $\mu$ is small. 
\medskip

A translation $T$ is made to center the coordinates at $P_2$  by moving $P_1$ to $e_1=(0,0,0)$ and $P_2$ to the origin:
$
q_1'=q_1- m_1, ~q_2'=q_3, ~q_3'=q_3,~ p_2'=p_2 - m_1, ~ p_1'=p_1, ~ p_3'=p_3.
$  
This is a symplectic map yielding a Hamiltonian,
$H(T^{-1}(q', p'))$. We add a constant, $\frac{(1-\mu)^2}{2}$ to this Hamiltonian as this does not change the Hamiltonian vector field, and end up with $H_m(q', p')= H(T^{-1}(q', p')) + \frac{(1-\mu)^2}{2}$.  For simplicity of notation, we replace $q', p'$ by $q, p$, not to be confused with previous notation.

Thus, the spatial lunar problem in translated $P_2$ centered coordinates can be represented as a Hamiltonian system, with Hamiltonian, 
\begin{equation}
H_m =  \frac{1}{2}|p|^2-\frac{\mu}{|q|} + \omega(q_1p_2-q_2p_1)-(1-\mu)\bigg(\frac{1}{\sqrt{(q_1+1)^2+q_2^2+q_3^2}}+q_1\bigg).
\label{eq:Hm}
\end{equation} 
The flow is given by,
\begin{equation}
\dot{q} = {H_m}_p , \quad \dot{p} = - {H_m}_q  .
\label{eq:HmDE}
\end{equation}
The energy surface $\Sigma_H$ becomes,
\begin{equation}
\Sigma_{H_m}= \{ (q,p) \in \mathbb{R}^6 | H_m(q,p) = -h \}, 
\label{eq:Hm-Surface}
\end{equation}
In order to study the flow in the coordinates $(q,p)$ defined by $H_m$ near $P_2$, we magnify the flow near $P_2$ by the map, $M: (q,p) \rightarrow (\hat{q}, \hat{p})$,
\begin{equation}
M: \quad \hat{p} = \mu^{-{1\over{3}}}p, \quad  \hat{q} = \mu^{-{1\over{3}}}q.
\label{eq:Magnify}
\end{equation}
Thus, for $\mu$ small, as we are assuming for this paper, the coordinates $(\hat{q}, \hat{p})$ are defined in a magnified neighborhood about $P_2$. This implies that when solutions are found in these coordinates, in the original coordinates $(q,p)$, the solutions lie close to $P_2$ as determined by \eqref{eq:Magnify}.

The transformation given by \eqref{eq:Magnify} is not symplectic. In order to obtain a Hamiltonian system in the coordinates $(\hat{p},\hat{q})$, it is noted that $M$ is conformally symplectic with constant conformal factor $\mu^{2/3}$. It is verified that a new Hamiltonian incorporating this magnification is given by,
$$
H^\mu(\hat{q},\hat{p})=\mu^{-2/3}(H_m(M^{-1}(\hat{p},\hat{q}))+1-\mu).
$$  This can be simplified using a Taylor expansion. This follows since, 
\begin{equation}
\label{eq:lunar_ham_full}
\begin{split}
H^\mu(\hat{q},\hat{p})&=\frac{1}{2}|\hat{p}|^2-\frac{1}{|\hat{q}|}+ \omega(\hat{q}_1\hat{p}_2-\hat{q}_2\hat{p}_1)\\
&\phantom{=} -\frac{1-\mu}{\mu^{2/3}}\bigg(\frac{1}{\sqrt{(\mu^{1/3}\hat{q}_1+1)^2+\mu^{2/3}\hat{q}_2^2+\mu^{2/3}\hat{q}_3^2}}+\mu^{1/3}\hat{q}_1-1\bigg).
\end{split}
\end{equation}
It is verified that the term within the square root in the last term of $H^\mu$ can be written as $1+2\mu^{1/3}\hat{q}_1+\mu^{2/3}|\hat{q}|^2$. Setting $x = 2\mu^{1/3}\hat{q}_1+\mu^{2/3}|\hat{q}|^2$, and using the formula,
$$\frac{1}{\sqrt{1+x}}=1-\frac{x}{2}+\frac{3x^2}{8}+\mathcal{O}(x^3)$$
for $|x| < 1,$ which is satisfied for $\mu$ sufficiently small,
it is verified that,
\begin{equation}
\label{eq:Hmu}
\begin{split}
H^\mu(\hat{q},\hat{p}) &=
\frac{1}{2}|\hat{p}|^2-\frac{1}{|\hat{q}|}+ \omega(\hat{q}_1\hat{p}_2-\hat{q}_2\hat{p}_1) +\frac{1}{2}|\hat{q}|^2-\frac{3}{2}\hat{q}_1^2+\mathcal{O}(\mu),
\end{split}
\end{equation}
where the term $\mathcal{O}(\mu)$ is a real analytic function of $\hat{q}$.
The Hamiltonian flow is given by, 
\begin{equation}
\dot{\hat{q}} = {{H^\mu}_{\hat{p}} } , \quad \dot{\hat{p}} = - {H^\mu}_{\hat{q}}.
\label{eq:HmuDE}
\end{equation}
The Hamiltonian flow takes place on fixed energy surfaces, 
\begin{equation}
\Sigma_{H^{\mu}}(h) = \{ (\hat{q},\hat{p}) \in \mathbb{R}^6 | H^{\mu}(\hat{q},\hat{p}) = -h \}. 
\label{eq:Hmu-Surface}
\end{equation}
It is remarked that setting $\mu = 0$ defines Hill's Problem. For small $\mu$, the (rescaled) restricted three-body problem 
represents a perturbation of order $\mu^{1\over{3}}$.

The function $H^\mu(\hat{p},\hat{q})$ has the form of a Hamiltonian for a perturbed rotating Kepler problem similar to what occurs in the restricted three-body problem about the primary mass point. As was studied in \cite{Belbruno:1981a} for motion about the primary mass point $P_1$, this system for $\mu = 0$ has the $\hat{q}_3$-axis as an invariant submanifold for the flow for $\hat{q}(t), \hat{p}(t)$.  

Let $\hat{\phi}^*(t)$ represent the solution on the $\hat{q}_3$-axis for $\mu = 0$. Setting $\hat{q}_k = 0, \hat{p}_k = 0 , k= 1,2,$ one obtains an integrable Hamiltonian system of 1 degree of freedom in the variables $(\hat{q}_3,\hat{p}_3)$.
After regularizing collisions, $P_0$ moves to some finite distance $d(h)$ from the origin, where $\dot{\hat{q}}_3 = 0$, and then falls back to $P_2$ for another collision.
It continues to do this in a periodic fashion for all time. This defines a periodic consecutive collision orbit 
with energy $-h$ and with period $T^*$. As $h$ varies,  $T^* = T^*(h)$ varies in a continuous manner. In contrast to the family studied in \cite{Belbruno:1981a}, the family studied here exists for all energies $h$. In particular, as $h$ increases to $0$, the distance $d$ remains bounded. 
 The family extends to positive energy $h$, and $d$ tends to $\infty$ as $h$ goes to $\infty$. 
\footnote{It is noted that there are two consecutive collision orbits, one on the positive $\hat{q}_3$-axis and the other on the negative $\hat{q}_3$-axis. We just consider the orbit on the positive axis, without loss of generality. }
\medskip

\noindent
In summary, the family of consecutive collision orbits, for $\mu =0$ for System \eqref{eq:HmuDE} on the $\hat{q}_3$-axis with frequency $T^*(h)$ lies on the energy surface $\Sigma_{H^{\mu}}(h)|_{\mu =0}$. This family is denoted by $ \mathcal{F}^*(h)$ and an orbit of this family is labeled $\hat{\phi}^*(t,h)$. This orbit moves a maximal distance $d(h)$. (In the original system given by \eqref{eq:HmDE} one has a similar family of consecutive collision orbits for $\mu =0$ on the energy surface $\Sigma_{H_m}$ given by \eqref{eq:Hm-Surface} which move close to $P_2$ to within $\mathcal{O}(\mu^{1/3} )$.
\medskip\medskip\medskip

We will prove the following theorem,
\begin{theorem}
\label{thm:per_orbit_perturbed_hill}
On each fixed energy surface $\Sigma_{H^{\mu}}(h)$ for System \eqref{eq:HmuDE}, there exists a unique periodic orbit, $\hat{\phi}(t,h,\mu)$, for $\mu$ sufficiently small, where $\hat{\phi}(t,h, 0) = \hat{\phi}^*(t,h)$  and whose period $T(\mu)$ continuously tends to $T(0) = T^*$, provided the orbit $\hat{\phi}(t,h, 0)$ is non-degenerate. ($\Delta \neq 0$ (see \eqref{eq:Delta}.))
The orbits of this family, $ \mathcal{F}(h, \mu)$, are symmetric to the $\hat{q}_2\hat{q}_3$-plane, and $ \mathcal{F}(h, 0) = \mathcal{F}^*(h)$ . 
\end{theorem}
We will refer to the orbits of this theorem as \emph{polar orbits}.

\begin{remark}
We shall see in the numerical section that the non-degeneracy condition appears to fail only twice: once for $h\leq 0$, or once for $h>0$, see Figures~\ref{fig:eigenvalues_lunar} and \ref{fig:Bifurcation}.
This is in contrast to the polar orbit in the rotating Kepler problem, which becomes degenerate infinitely many times.
This may seem surprising given the similarities between the Hamiltonians of Hill's lunar problem and the rotating Kepler problem. However, we point out several important differences:
\begin{enumerate}
\item The rotating Kepler problem is completely integrable, whereas Hill's lunar problem is not: the additional terms in the potential describe a tidal and centrifugal force.
\item The period of the polar orbit in the rotating Kepler problem goes to infinity as the Jacobi energy goes to $0$. For Jacobi energy $h>0$, a regularized orbit moving on the $z$-axis escapes to infinity.
In Hill's lunar problem, the region consisting of the intersection of the Hill's region with the $z$-axis (containing the origin) is bounded for all energies $h$.
As a result consecutive collision orbits in Hill's lunar problem are periodic for all energies $h$.

\item  From a physical point of view there is a large difference between these problems. In the rotating Kepler problem, the rotational term $q_1p_2-q_2p_1$ is  due to the rotating coordinate system centered at the larger primary in $0$.
In Hill's lunar problem, the center of rotation is infinitely far away: the physical meaning of the rotational term is hence more complicated, resulting in additional terms corresponding to a tidal/centrifugal force.

\item The orbits in Hill's lunar problem become unstable for large $h$. Intuitively, the instability in the polar orbits in Hill's lunar problem for large energies is easy to understand: for sufficiently large energy $h$ they spend a considerable time away from the smaller primary centered at the origin, so that the tidal forces can destabilize the orbit.
\end{enumerate}

\end{remark}

\subsubsection*{Distance of orbits to the Moon}
The map $M$ given by \eqref{eq:Magnify} scales the coordinates of the periodic orbits by $\mu^{{1\over{3}}}$ when mapping back to the original $(q,p)$ coordinates of \eqref{eq:HmDE}. Thus for $\mu$ small, the periodic orbits remain close to $P_2$ to within the distance, 
$\rho = \mathcal{O}(\mu^{{1\over{3}}})$. This distance, however, is significant and can extend to the $L_1, L_2$ Lagrange points. 
\subsubsection*{Bounded period and existence for all energies}
We estimate the period of the polar orbit in Hill's lunar problem. The rotational term drops out on the $z$-axis, so the energy equals
$$
E_H=\frac{1}{2}\dot z^2+V(z),
$$
where we take Hill's lunar potential restricted to the $z$-axis, given by $V(z)=-\frac{1}{z} + \frac{1}{2}z^2$.
Fix the energy to the value $h$.
The particle moves between $z=0$ and $z=d(h)$, where $d(h)$ is a solution to $V(z)=h$. This equation is equivalent to the cubic equation
$$
z^3-2hz-2=0,
$$
which clearly has a unique, positive solution, which can be found with Cardano's formula.
Using the energy, we can compute the speed, and find for the period
$$
T^*=2\int_0^{d(h)}\frac{1}{\sqrt{2h-V(z)}}\mathrm{d}z
=2\int_0^{d(h)}\frac{1}{\sqrt{2h+\frac{2}{z}-z^2}}\mathrm{d}z.
$$
This integral can be evaluated exactly using elliptic integrals as one may verify with a computer algebra system such as \texttt{Maple}.
The expression is not too illuminating, and we will only establish a period bound here.
We compute
\[
\begin{split}
T^* & \underset{u=z/d(h)}{=} 2\int_0^1 \frac{1}{\sqrt{2c+\frac{2}{d(h) u} - d(h)^2 u^2}} d(h) \mathrm{d}u 
\leq 2\int_0^1 \frac{1}{\sqrt{\frac{2c}{d(h)^2}+\frac{2}{d(h)^3}-u^2  } } \mathrm{d}u\\
&=
2\int_0^1 \frac{1}{\sqrt{1-u^2}} \mathrm{d}u=\pi
\end{split}
\]
In other words, the polar orbit in Hill's lunar problem has a uniform period bound holding for all $h$. Furthermore, this period bound is so small that the polar orbit always closes up before the rotational term can finish even one revolution. We shall see that the orbit still becomes degenerate due to the tidal and centrifugal force terms.

\section{Proof}
\label{sec:Proof} 

In this section we prove Theorem 1.  It is necessary to regularize the flow since the consecutive collision orbits for $\mu = 0$ collide with $P_2$. After that is done, making use of the symmetry of the orbits, allows a section to be defined. The existence of the periodic orbit family then results from an application of the implicit function theorem.

The main Hamiltonian, $H^\mu(\hat{q},\hat{p})$, \eqref{eq:Hmu}, in this paper that we want to regularize has a form similar to the Hamiltonian in \cite{Belbruno:1981a}(Equ. 1, p. 397). The differences are that roles of the $q_1$ and $q_2$ axis are reversed, and \eqref{eq:Hmu} has the extra term, 
$$E = +\frac{1}{2}|\hat{q}|^2-\frac{3}{2}\hat{q}_1^2 .
$$
This term is smooth at $\hat{q} = 0$, and will not affect the regularization. We will use the results of \cite{Belbruno:1981a} often, and henceforth refer to it as B81 for the convenience of the reader.

The idea of regularization is to make a symplectic transformation of the coordinates,
$$ \hat{q} = F(P,Q), \hspace{.1in} \hat{p} = G(P,Q),$$ 
where $Q \in \mathbb{R}^3, P \in \mathbb{R}^3$ and a transformation of time $t \rightarrow s$, so that in the new coordinates, the Hamiltonian flow is well defined at collision.

In three degrees of freedom, regularizations are considerably more complex than in two degrees of freedom, where, for example, the Levi-Civita transformation can be readily applied. A regularization for three degrees of freedom
is developed and applied in B81 that is ideal for the collision at $P_2$, since, as noted, the Hamiltonian $H$ in B81 is very close to $H^\mu(\hat{q},\hat{p})$.  

The regularization in B81 is a higher dimensional M\"obius transformation. It is represented in a clear manner by 
defining a Jordan algebra that serves as a generalization of the complex numbers. This is a nonassociative algebra defined on the space $\mathbb{A}_{n}$ which is isomorphic to $\mathbb{R}^{n+1}$ as a vector space. Its product structure is defined as follows. Write $z \in  \mathbb{A}_{n}$ as 
$$z = z_0 + i_1z_1 +i_2z_2 + ... + i_n z_n,$$ 
$z_i \in \mathbb{R}^1$ for $i = 1,2,...,n$. The $\mathbb{R}$-linear multiplication is then defined by imposing 
$$i_\alpha i_\beta = - \delta_{\alpha\beta}.$$  
Conjugation is defined as 
$$
\bar{z} = z_0 - i_1z_1 -i_2z_2 - ... - i_n z_n.
$$ 
One then obtains, $z\bar{z} = |z|^2 = z_0^2 + z_1^2 + ... + z_n^2$. Division is defined as 
$$
\frac{1}{z} = \frac{\bar z}{|z|^2}
$$
Although this algebra is commutative, it is non-associative. A measure of the non-associativity is the 'associator' $a = x(yz) - (xy)z$, $x,y,z$ are each in $\mathbb{A}_{n}$.  Since our variables are of three components, the case $n=2$ is of relevance, where $a = (x_2z_1 - x_1z_2)(i_1 y_2 - i_2 y_1)$.  Further details are in B81. 

With coordinates in this Jordan algebra $\mathbb{A}_{2}$ we will obtain a simple form for a symplectic transformation that regularizes collisions. Since the axes in this paper differ from those in B81, we first interchange the first two components of $\vec q, \vec p$: $\vec q=(\hat q_2,\hat q_1,\hat q_3)$ and $\vec p=(\hat p_2,\hat p_1,\hat p_3)$. After that, we use the transformation
$$
\vec{p} = {{1+ P}\over{1 -P}}, \hspace{.1in} \vec{q} = {{Q}\over{2}}(1- \bar{P})^2 - (QP)\bar{P} + Q(P\bar{P}).
$$
As in B81, we will see that the Hamiltonian $\Gamma$ given by
$$ 
\Gamma = (1/2)|P - (1,0,0)|^2|Q| (H^\mu +h) 
$$
regularizes the energy level $H^\mu=-\mu$; we use the time transformation $t \rightarrow s$ given by $ t = \int^s|\vec{q}| d\tau = (1/2)\int^s|P - (1,0,0)|^2|Q|d\tau $ for the rescaling of the Hamiltonian. The new level set of interest, corresponding to $H^\mu=-\mu$ will be denoted by
$$
\Sigma_\Gamma = \{(P,Q) \in \mathbb{R}^6| \Gamma = 0 \}.
$$

The new Hamiltonian has the form of B81(Equ. 22) with the addition of the term $\tilde{E}(Q,P)$. Here $\tilde{E}(Q,P)= \frac{1}{8} |P - (1,0,0)|^4|Q|^2 - \frac{3}{2}f^2(Q,P)$, where $f(Q,P)$ is the first component of the transformation of $ \vec{q}$.
\[
\begin{split}
\Gamma &=  {|Q|\over{4}}\{|P+(1,0,0)|^2 + 2(h + \omega \alpha(Q,P))|P-(1,0,0)|^2 \} \\
&\phantom{=} - 1 +
 (1/2) \tilde{E}(Q,P)|P - (1,0,0)|^2|Q| + \mathcal{O}(\mu). 
\end{split}
\] 
The Hamiltonian flow is defined by, 
$$
Q^{'} = \Gamma_P , \hspace{.1in} P^{'} = - \Gamma_Q , 
$$
where $\phantom{ }^{'} \equiv {d\over{ds}} $. 
To check that the Hamiltonian flow $X(s) = (Q(s), P(s))$ is regular at collision orbits, we note that collision occurs when $ \hat{\phi}^*(t,h) = (\hat{q}^*(t),\hat{p}^*(t))$ tends to $(0,0,0;0,0, \infty)$. In the coordinates $Q,P$, any collision point corresponds to $P=(1,0,0)$ and $|Q|=1$; the collision point $X(0)=(0,0,-1;1,0,0)$. We see that the Hamiltonian $\Gamma$ is smooth near collision points, so the collision orbit, which we label by $X^*(s)$, is indeed regularized and becomes an honest periodic orbit on $\Sigma_\Gamma$. We will denote its period by $S^*$.

The existence of a unique periodic orbit $X(s,X_0, \mu)$ near $X^*(s)$, $X(0,X_0, \mu) = X_0$, of period $S$ near $S^*$ for $\mu$ sufficiently small, is obtained by the implicit function theorem, applied to the subset of symmetric orbits, as we shall now see.

The Hamiltonian flow is symmetric with respect to the $\hat q_1,\hat q_3$-plane, because the Hamiltonian $H^\mu$ is invariant under the map 
$$ 
t \rightarrow -t, \hspace{.1in}  (\hat{q}_1,\hat{q}_2,\hat{q}_3,\hat{p}_1,\hat{p}_2, \hat{p}_3) \rightarrow (\hat{q}_1,-\hat{q}_2,\hat{q}_3,-\hat{p}_1,\hat{p}_2, -\hat{p}_3).
$$ 
Keeping in mind our interchange of components, this implies that the involution
$$ 
s \rightarrow -s, \hspace{.1in}  (Q_1,Q_2,Q_3,P_1,P_2, P_3) \rightarrow (-Q_1,Q_2,Q_3,P_1,-P_2, -P_3).
$$
preserves $\Gamma$. Solutions that are symmetric with respect to this involution are characterized by the condition,
$$
Q_1(0) = P_2(0) = P_3(0) = 0.
$$
This means that symmetric solutions are characterized by three initial values 
$Q_2(0), Q_3(0), P_1(0).$ 
This can be reduced to two initial values on the energy surface near $X^*(0)$. Namely at $X^*(0)$ we can verify that 
\begin{equation}
  { {\partial \Gamma}\over{\partial Q_3 } }= -1\neq 0. 
	\label{eq:NonZeroDeriv}
	\end{equation}
\eqref{eq:NonZeroDeriv} also proves that the flow near $X^*(0)$ is transverse to $X^*(0)$. 
Thus, by the implicit function theorem, near $X^*(0)$, we can eliminate the $Q_3$ coordinate and characterize symmetric solutions by two initial values, 
$$Q_2(0), P_1(0).$$

\noindent
In addition, the time of intersection of solutions near $X^*(0)$ can be determined from \eqref{eq:NonZeroDeriv}  as $s = s(Q_2, Q_3, P_1)$ by the implicit function theorem for $\mu$ sufficiently small. 

This defines a three-dimensional surface of section 
$$\mathcal{S} = \{ Q_2, Q_3, P_1 | Q_1(0) = P_2(0) = P_3(0) = 0 \} .$$

A solution starting on this section at time $s =0$, then reintersecting the section at time $s = S/2$ yields a symmetric periodic orbit of period $S$. 
This is satisfied by the consecutive collision orbit, $X^*(s)$ with $s = S^*/2$. For this to be satisfied
near  $X^*(s)$  for small $\mu$ by a solution $X(s,X_0, \mu)$ yields a periodic orbit of period $S(\mu)$ near $S^*$, such that $S(0) = S^*$. This can be satisfied provided the determinant 
\medskip

\begin{equation}
\Delta = det  { {\partial(Q_1, P_2, P_3)}\over {\partial(S, Q_2(0), P_1(0)) } } 
\label{eq:Delta}
\end{equation}
\medskip

\noindent does not vanish at $\mu = 0, S = S^*/2, X(0) = X^*(0)$. When 
$\Delta \neq 0$ 
is satisfied, we say that the periodic orbit $\hat{\phi}(t,h,0)$ is {\it non-degenerate}. It is numerically shown next, in Section \ref{sec:Numerical}, that $\hat{\phi}(t,h,0)$ is non-degenerate except for two values of $h$.
\medskip

\noindent

\medskip

\noindent
This concludes the proof of Theorem 1.


\section{Numerical Results}
\label{sec:Numerical}
We start with a summary of the numerical results. Throughout, we will be comparing Hill's lunar problem with the rotating Kepler problem. The reason for this is twofold.
\begin{itemize}
\item The same type of polar orbit has been studied before in B81. Comparison will hence clarify differences and similarities.
\item We can consider the two types of polar orbits as part of a larger family of periodic orbits in the restricted three-body problem. Both the rotating Kepler problem and Hill's lunar problem are limit cases where the polar orbit has a particularly simple form.
\end{itemize}
Here is a list of the main results. We will write $H$ for Hill's lunar problem and $K$ for the rotating Kepler problem.
\subsubsection{Stability properties for fixed $\mu$:}
\begin{enumerate}
\item[H:] The polar orbit goes through four bifurcations for $h\in (-\infty,\infty)$: they are a period doubling bifurcation, a simple degeneracy, another simple degeneracy, and a period halving bifurcation.
The polar orbit is elliptic for $h<-1.03$ and complex hyperbolic for $h>0.11$.\footnote{A quick overview of the terminology: by elliptic we mean two conjugate eigenvalues on the unit circle. This implies a weak form of stability: nearby orbits cannot escape quickly.

By hyperbolic we mean two real eigenvalues: $\lambda$ and $1/\lambda$. We add ``negative'' to indicate that $\lambda<0$.
The return map in the spatial problem has four eigenvalues, satisfying the symmetry property: if $\lambda$ is an eigenvalue, then so are $\bar \lambda$, $1/\lambda$ and $1/\bar \lambda$.
This leaves an one additional cases in this dimension, namely none of the eigenvalues are purely real, nor do they lie on the unit circle: we will call this complex hyperbolic.

All forms of hyperbolicity implies instability in the sense that nearby orbits tend to escape quickly: how quickly depends on the absolute value of the largest eigenvalue.
}
\item[K:] The polar orbit goes through infinitely many simple degeneracies for $h \in (-\infty,0)$, and the orbit is elliptic for all $h<0$. Simple degeneracies occur whenever the period of the polar orbit is a multiple of $2\pi$, the rotation period of the coordinate system.
\end{enumerate}

\subsubsection{Variation of the family $\mathcal F(h, \mu)$ with fixed $\mu$ and varying $h$}
\label{sec:varying_energy}
We consider small deformations of Hill's lunar Hamiltonian, i.e.~small $\mu$ in $H^\mu$ and of the rotating Kepler problem, i.e.~$\mu$ close to $1$ in $H_m$.
\begin{enumerate}
\item[H:] The polar orbit starts out as a very eccentric ellipse, staying near the $z$-axis closely: the projection to the $xy$-plane looks like an oval.
After becoming degenerate twice, the orbit starts to develop a cusp in the $yz$-projection. 
\item[K:] The orbit also starts out as a very eccentric ellipse, hugging the $z$-axis: the projection to the $xy$-plane looks like an oval for very negative energy. As $h$ increases, the shape ceases to be convex and the orbit becomes degenerate.
With each degeneracy, the winding number of the polar orbit around $0$ increases; in other words, the orbit accumulates loops as the Jacobi energy increases.
\end{enumerate}

\subsubsection{A bridge between polar orbits in the rotating Kepler problem and the restricted three-body problem near the light primary}
For energy $h\leq -1.50$ \footnote{We remind the reader that the critical energy for $\mu=0.5$ equals $-2.0$.} there is a bridge with constant Jacobi energy $h$ connecting polar orbits near the smaller primary, meaning small $\mu$ in $H_m$, to polar orbits in the rotating Kepler problem, meaning $\mu=1$ for the Hamiltonian $H_m$.
The orbits near the smaller primary can be continued to Hill's lunar problem after rescaling the coordinates.

For very negative Jacobi energy this bridge does not involve any dynamical transitions. For larger Jacobi energy, orbits near the light primary are hyperbolic, whereas orbits in the rotating Kepler problem are elliptic, so this bridge necessarily involves bifurcations.
We discuss the theory behind the existence of a bridge in Section~\ref{sec:Geometry}.

\begin{figure}[!htp]
  \centering
  \includegraphics[width=0.85\textwidth,clip]{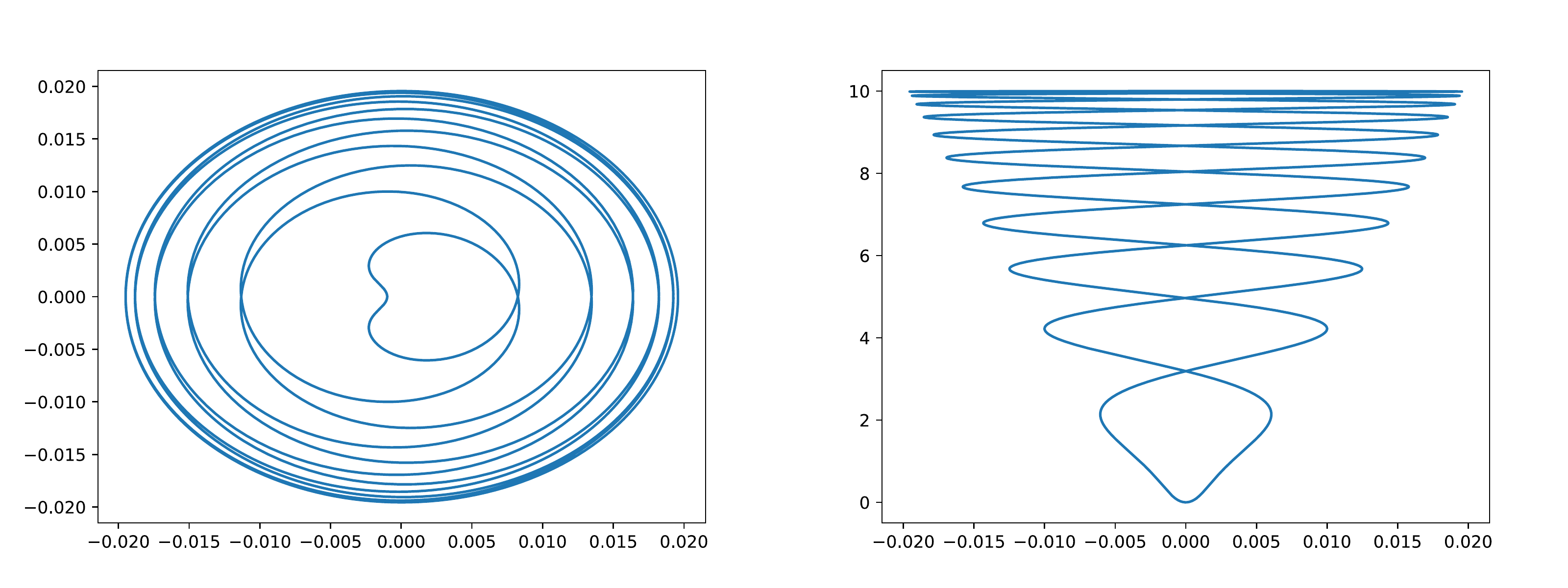}
  \caption{The $xy$- and $yz$-projection of an orbit in a partial bridge from the rotating Kepler: $\mu=0.999$, $h=-0.1$.}
  \label{fig:bridge4}
\end{figure}

\begin{figure}[!htp]
  \centering
  \includegraphics[width=0.85\textwidth,clip]{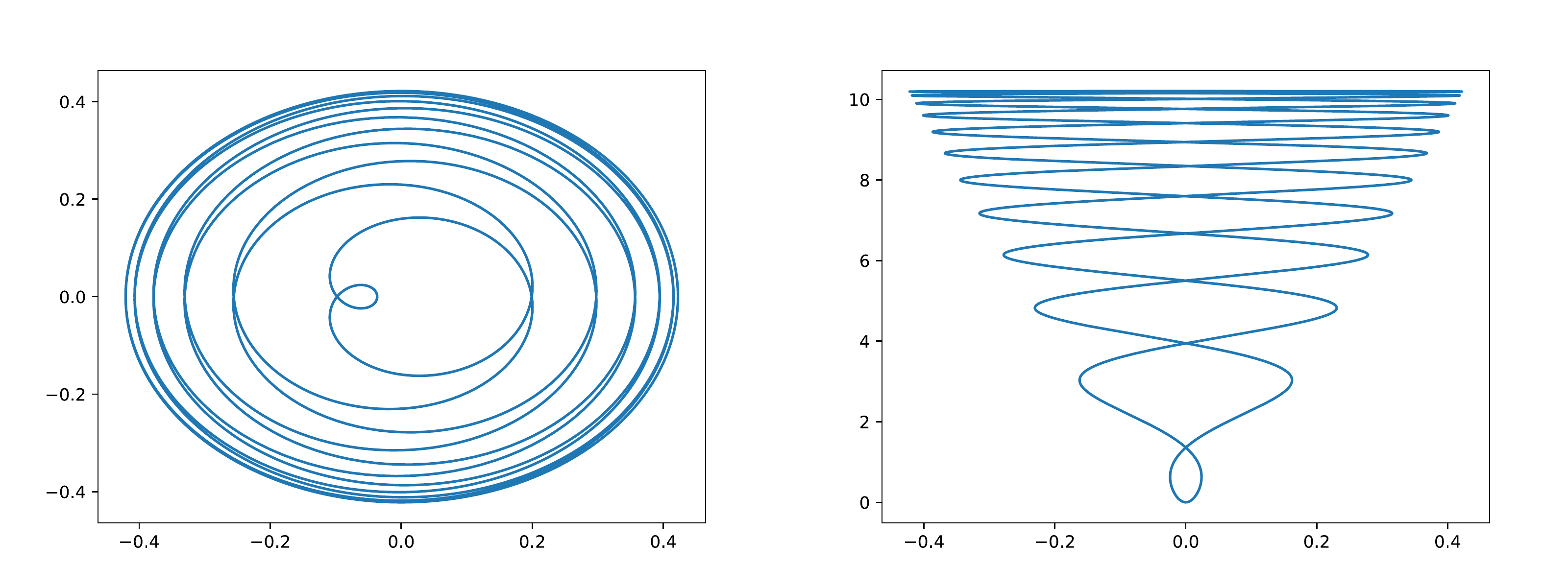}
  \caption{The $xy$- and $yz$-projection of an orbit in a partial bridge from the rotating Kepler: $\mu=0.963$, $h=-0.1$.}
  \label{fig:bridge5}
\end{figure}

Depending on the Jacobi energy, the bridge can also involve loops. A sample is given in Figures~\ref{fig:bridge4} and \ref{fig:bridge5}. Note that the loops only appear at higher Jacobi energies and bridges at such energies do not reach all the way to the lunar problem.

\subsubsection{Non-collision polar orbits in the Moon-Earth system}
\begin{figure}[!htp]
  \centering
  \includegraphics[width=0.85\textwidth,clip]{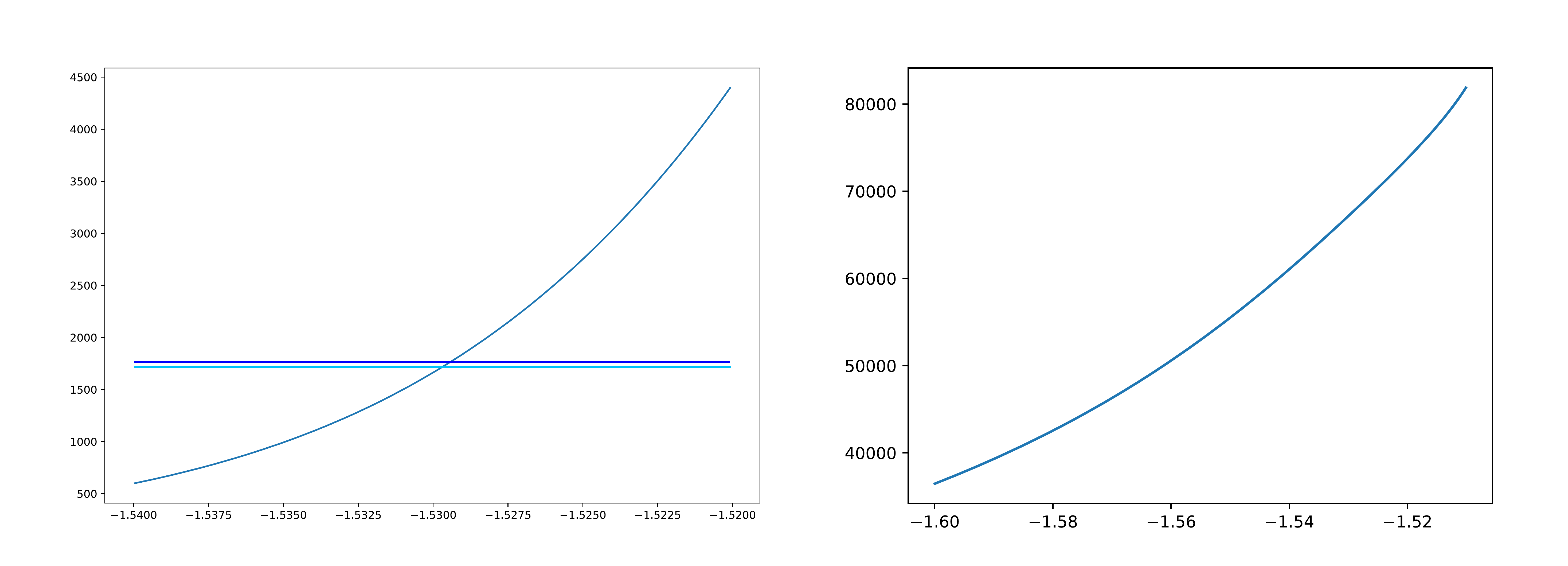}
  \caption{Periapsis (left) and apoapsis (right) in km as function of the Jacobi energy. At the light blue line, the orbit just hits the surface of the Moon (taken to have a radius of 1716 km), and at the dark blue line, the orbit reaches a periapsis that is least 50 km above the surface.}
  \label{fig:periapsis_apoapsis}
\end{figure}

We continue the polar orbit into the Moon-Earth system where the polar orbit turns out to be a physical non-collision orbit for sufficiently high Jacobi energy. A plot of the periapsis and apoapsis as function of the energy is given in Figure~\ref{fig:periapsis_apoapsis}. It also undergoes a period doubling bifurcation, making a transition from stable to unstable in the same energy range.

\subsection{Details concerning the numerical procedure}

\subsubsection{Regularization scheme }
We will use the Moser-Belbruno-Osipov regularization scheme to regularize the flow. We use the incarnation due to Moser, which we will refer to as just Moser regularization, and the incarnation due to Belbruno specialized to collision orbits as in B81, which we will refer to as Belbruno transform; this scheme was described in Section~\ref{sec:Proof}. Both schemes are detailed in the appendices.

As a short summary, the Moser scheme regularizes the energy hypersurface below the critical value to the unit cotangent bundle of the three-sphere and has the advantage that it is global, i.e.~no local charts are needed. However, to do computations in the Moser scheme we need to impose constraints to stay on this space, which we view as a submanifold of $T^*\mathbb{R}^{4}$; this leads to a slightly larger computational effort. 

The Belbruno transform uses a generalized M\"obius transformation, based on the Jordan algebra described in Section~\ref{sec:Proof}. For this regularization scheme, the advantages and disadvantages as described above are reversed. The scheme is local, and gives a chart, which some orbits could escape from. On the other hand, the Belbruno regularization does not require constraints.

\subsubsection{Integration scheme}
For numerical integration we have used a Taylor integrator with both variable stepsize and order. The typical order with a \texttt{double} and \texttt{long double}, which corresponds to about 15-16 digits and 18-19 digits precision, respectively, was between 20 and 30.
We have used three different implementations of the Taylor integrator: the Taylor translator described in~\cite{Jorba:Taylor}, the CAPD-libary, \cite{CAPD}, and a homegrown Taylor library.

To find periodic orbits for $\mu>0$, we made use of a local surface of section and the familiar homotopy method to follow solutions from $\mu=0$ to the desired value of $\mu$ in sufficiently steps of $\mu$.
We choose a linear surface of section perpendicular to the $z$-axis.
This is useful to follow the orbits for large parameter changes.
As usual, we followed the orbit until it crossed the surface of section, and found a more accurate intersection by normalizing the flow.

For the stability analysis of the polar orbit, we choose a symmetric surface of section in line with the proof in Section~\ref{sec:Proof}. This has the advantage that the symmetry properties can be exploited more effectively.

A final remark concerning the Hamiltonian: for the lunar problem we use the regularization of $H^\mu$, the Hamiltonian given in~\eqref{eq:lunar_ham_full}, but some care has to be taken to deal with the catastrophic cancellation in the final term.

\subsection{Detailed results}

\subsubsection{Non-degeneracy of the lunar and Kepler orbit for $h\leq 0$}
\label{sec:non-degeneracy}
To apply the theorem, we need to know whether the non-degeneracy condition holds. To check this, we numerically compute the linearized return map transverse to the flow. We can represent this as a symplectic $4\times 4$-matrix, so its eigenvalues have some symmetry properties. Namely, if $\lambda$ is an eigenvalue of a symplectic $4\times 4$-matrix, then $\bar \lambda$, $1/\lambda$ and $1/\bar \lambda$ are also eigenvalues (possibly equal). In general, this leaves a lot of possibilities. However, it turns out that the linearized return map is elliptic, i.e.~all eigenvalues lie on the unit circle, for very negative energies $h$.
The behavior for the lunar problem, which is of primary interest here, turns out to differ from the behavior in the rotating Kepler problem studied in B81. 

\begin{figure}[htp]
  \centering
  \includegraphics[width=0.65\textwidth,clip]{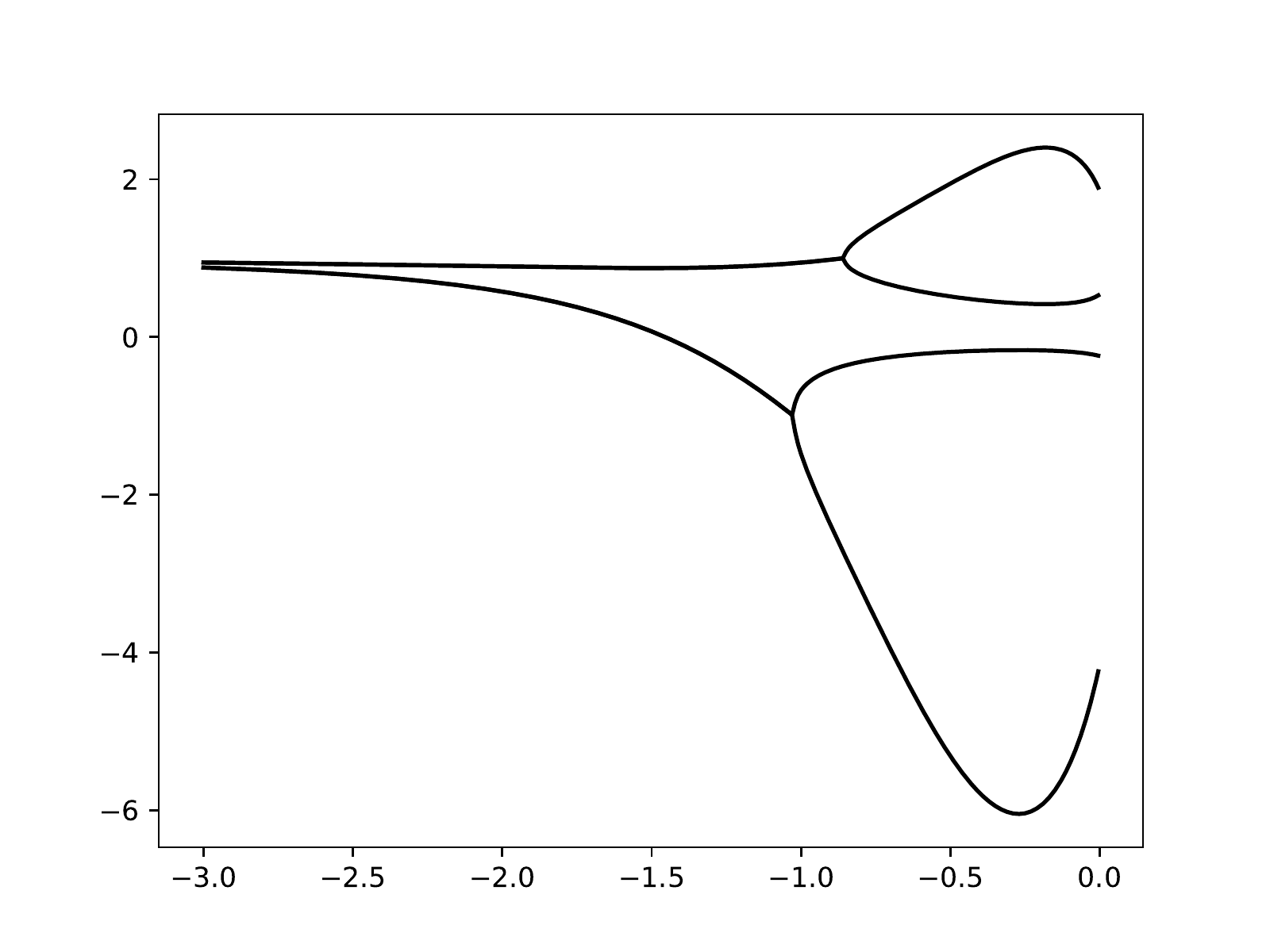}
  \caption{The real part of the four eigenvalues for $\mu=0$, the lunar problem, as a function of the energy $h$ in the Hamiltonian $H^{\mu}$ for $h\leq 0$. Bifurcations from an elliptic to a hyperbolic eigenvalue appear twice.}
  \label{fig:eigenvalues_lunar}
\end{figure}

See Figure~\ref{fig:eigenvalues_lunar} for the eigenvalues in the lunar problem. For very negative energies, the return map is elliptic. For $h\sim -1.03$, the orbit goes through a period doubling/halving bifurcation: the orbit goes from being purely elliptic to mixed elliptic/negative hyperbolic without becoming degenerate; its double cover does become degenerate.

At energy $h\sim -0.86$, the orbit itself becomes degenerate resulting in a positive hyperbolic/negative hyperbolic pair of eigenvalues.

\begin{figure}[!htp]
  \centering
  \includegraphics[width=0.6\textwidth,clip]{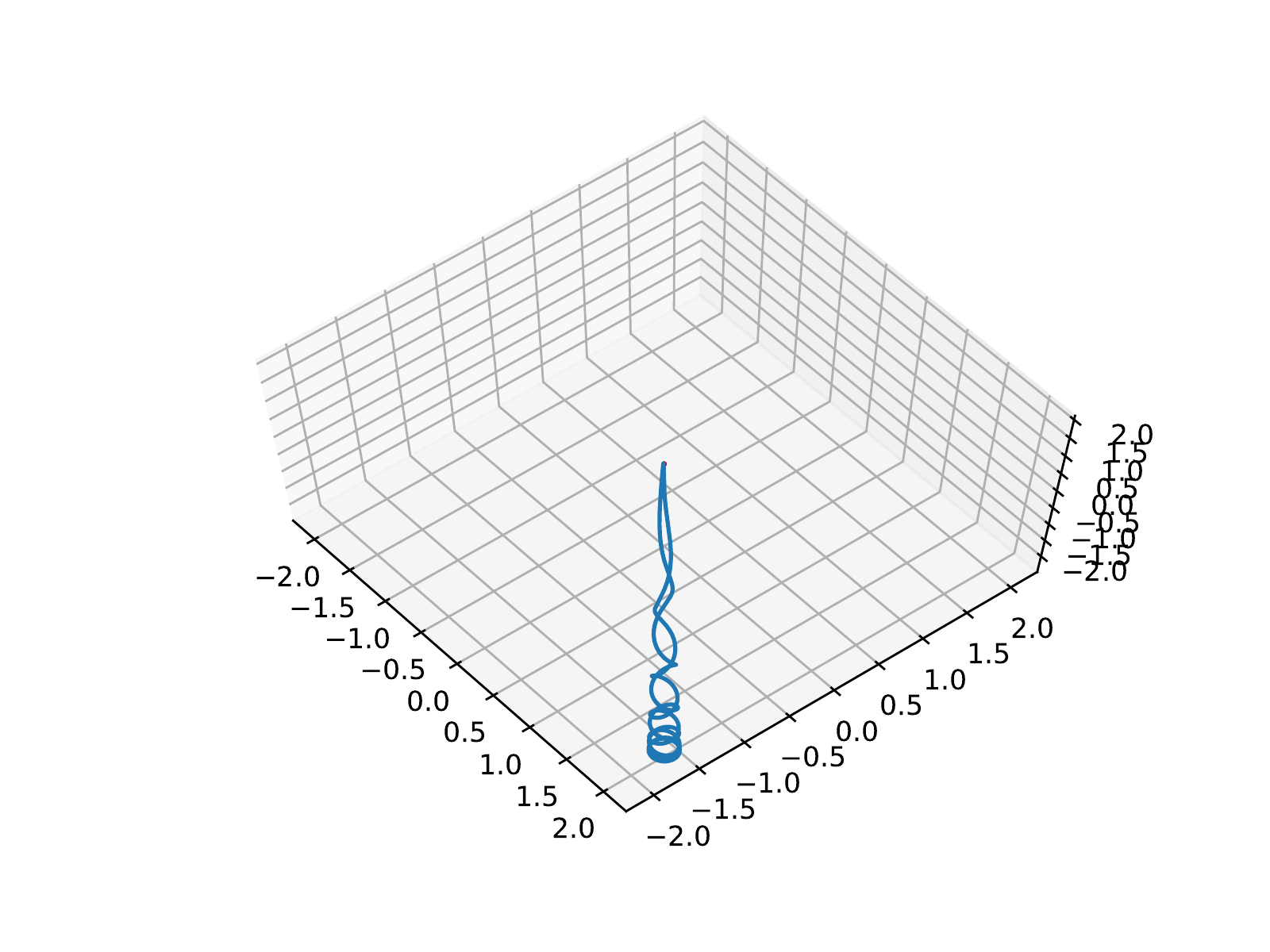}
  \caption{A periodic polar orbit in RTBP for large mass ratio, close to the rotating Kepler problem. This orbit has picked up many loops as described in Section~\ref{sec:varying_energy}}
  \label{fig:multi_loop}
\end{figure}

This is in contrast to the situation for the rotating Kepler problem, where the orbit stays elliptic up to $h=0$.
Furthermore, the polar orbit in the rotating Kepler problem becomes degenerate infinitely many times as the energy goes to $0$, and its behavior changes every time when it does. This results in loops appearing in the perturbations of the rotating Kepler problem. This behavior was found by Belbruno in B81. We include an illustration for the convenience of the reader in Figure~\ref{fig:multi_loop}.
Most of the other illustrations of the orbits will involve only projections to the $xy$- and $yz$-plane.
\begin{figure}[htp]
  \centering
  \includegraphics[width=0.65\textwidth,clip]{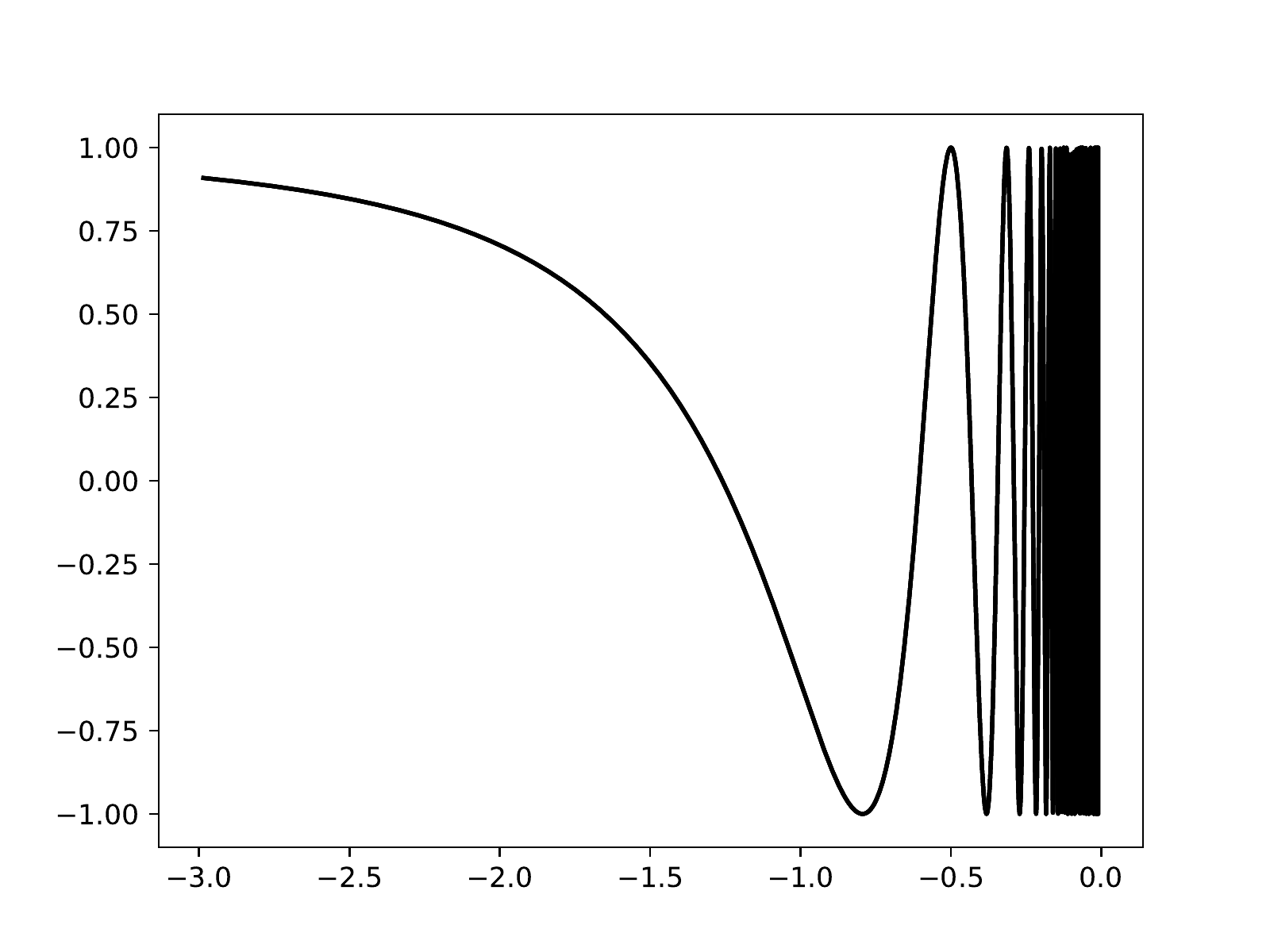}
  \caption{The real part of two eigenvalues for $\mu=1$, the rotating Kepler problem, as a function of the energy $h$ in the Hamiltonian $H_m$. All eigenvalues are elliptic, and bifurcations occur whenever an eigenvalue passes through 1. }
  \label{fig:eigenvalues_kepler}
\end{figure}
\begin{remark}
The plot in Figure~\ref{fig:eigenvalues_kepler} was obtained through numerical means, and we want to point out that one can obtain an analytical expression for the eigenvalues of the linearized return map.
\end{remark}

\subsubsection{Stability properties of the polar orbit in the lunar problem for $h>0$}
\begin{figure}[!htp]
  \centering
  \includegraphics[width=0.25\textwidth,clip]{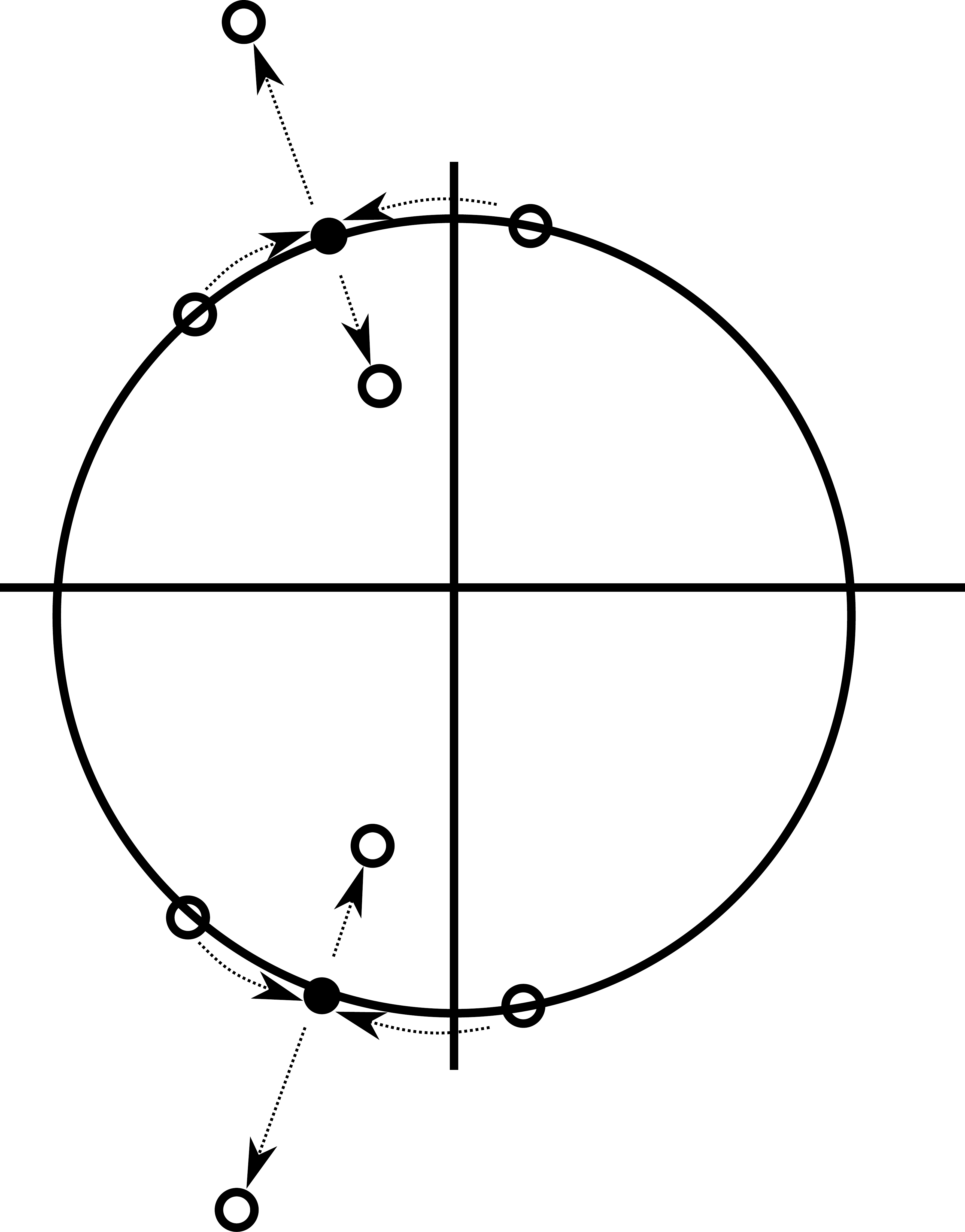}
  \caption{Collision of eigenvalues}
  \label{fig:eigenvalues_collision}
\end{figure}

In Hill's lunar problem, the polar orbit remains a periodic orbit for $h>0$, and its stability properties there are very interesting. 
To understand the situation, recall that eigenvalues of a symplectic matrix come with symmetries as mentioned in Section~\ref{sec:non-degeneracy}.
In the spatial problem, there are four eigenvalues, and an orbit can lose stability without becoming degenerate by the following mechanism:
\begin{enumerate}
\item at parameter $h_0$ all eigenvalues are elliptic, i.e.~on the unit circle.
\item as the parameter $h\to h_1$ the eigenvalues stay elliptic, but they collide in pairs: i.e. there are only two distinct eigenvalues.
\item for $h>h_1$ the eigenvalues move off the unit circle as sketched in Figure~\ref{fig:eigenvalues_collision}.
\end{enumerate}

It turns out that this mechanism occurs in for Hill's lunar problem for $h>0$.
\begin{figure}[!htp]
  \centering
  \includegraphics[width=0.75\textwidth,clip]{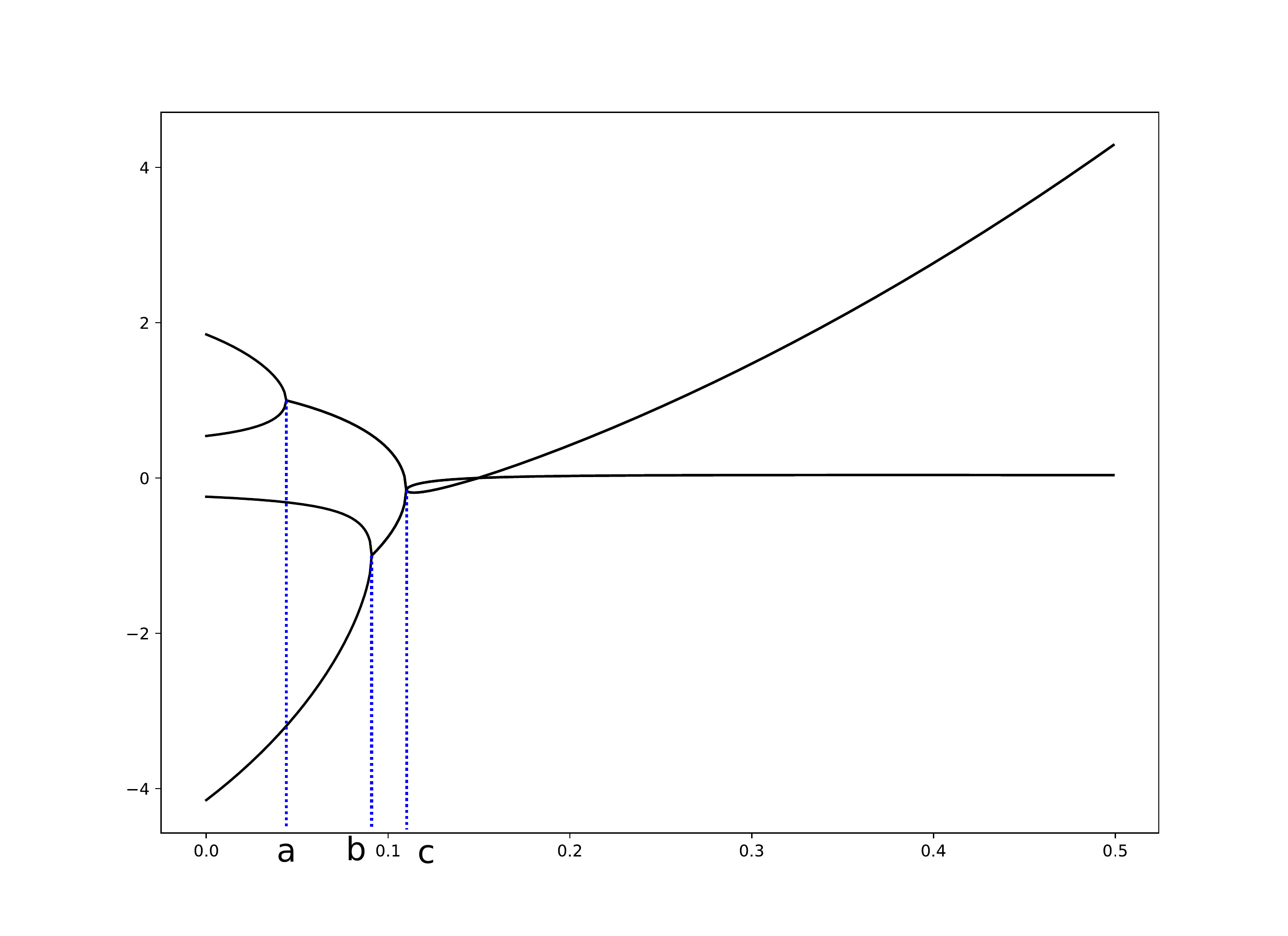}
  \caption{Bifurcation of the polar orbit for $h>0$ in Hill's lunar problem: the real part of the four eigenvalues of the linearized return map.}
  \label{fig:Bifurcation}
\end{figure}

We briefly explain the bifurcation points in Figure~\ref{fig:Bifurcation}
\begin{enumerate}
\item[a] the orbit becomes degenerate and goes from being hyperbolic/negative hyperbolic to elliptic/negative hyperbolic for $h\sim 0.044$.
\item[b]  the orbit goes through a period doubling/halving bifurcation: it goes from elliptic/negative hyperbolic to elliptic/elliptic for $h\sim 0.091$.
\item[c] the orbit goes through an eigenvalue collision for $h\sim 0.11$ and the eigenvalues move from away from the unit circle. They do not appear to return to the unit circle, so stability seems to be lost for large $h>0$.
Beyond point c, the eigenvalues move ``freely'' in the complex plane, so the  real part has then little meaning by itself. In particular, the additional ``intersection'' is not an intersection in the complex plane. 
\end{enumerate}

\begin{remark}
The same mechanism of losing stability takes place for small $\mu>0$.
\end{remark}
We summarize the precise results in the following proposition.
\begin{proposition}
Let $\gamma_h$ denote the polar orbit in Hill's lunar problem as a function of the energy $h$.
Then $\gamma_h$ undergoes the following bifurcations:
\begin{enumerate}
\item a period doubling/halving bifurcation for $h\in[-1.025245,-1.025225]$ and in the interval $[0.0909615,0.0909616]$.
\item a simple degeneracy for $h\in[-0.85556,-0.85555]$ and for $h\in[0.043843, 0.043844]$. 
\item an eigenvalue collision for $h\in[0.109989, 0.109990]$.
\end{enumerate}
\end{proposition}

\textbf{Proof:}
The argument is computer-assisted, and uses interval arithmetic to obtain rigorous error bounds.
The idea is to make the proof from Section~\ref{sec:Proof} more quantitative, and to use the linearized flow to obtain a tight enclosure of the return map.
By starting on the symmetric surface of section, we only need to compute half the periodic orbit. 

After obtaining an enclosure for the return map using the linearized flow, we compute the coefficients of the characteristic polynomial of the restriction of the linearized flow to a transverse slice. Let us denote this restriction by $\psi$. We have
$$
\chi(\psi)(x)=x^4-s_1(\psi)x+s_2(\psi)x^2-s_3(\psi)x+\det(\psi).
$$
Here $s_i(\psi)$ denotes the elementary symmetric polynomial in the roots of the characteristic polynomial of degree $i$, so $s_1(\psi)=Tr(\psi)$.
A $4\times 4$ symplectic matrix satisfies $\det(\psi)=1$ and $s_3(\psi)=s_1(\psi)$.
Hence we can compute the entire characteristic polynomial by just computing $s_1$ and $s_2$. Using the standard formula for a quadratic polynomial we find all roots with good error bounds (direct computation of the determinant gives far worse enclosures).
\hfill $\square$
\\

Based on numerical experiments, we obtain the following.
\begin{observation}
The orbit is elliptic for $h<-1.03$ and complex hyperbolic for $h>0.11$
\end{observation}
Although we have checked the statement for some finite intervals using interval arithmetic, we do not have a full proof that works for all energies.

\subsubsection{Following the orbits for varying energy}
\begin{figure}[!htp]
  \centering
  \includegraphics[width=0.85\textwidth,clip]{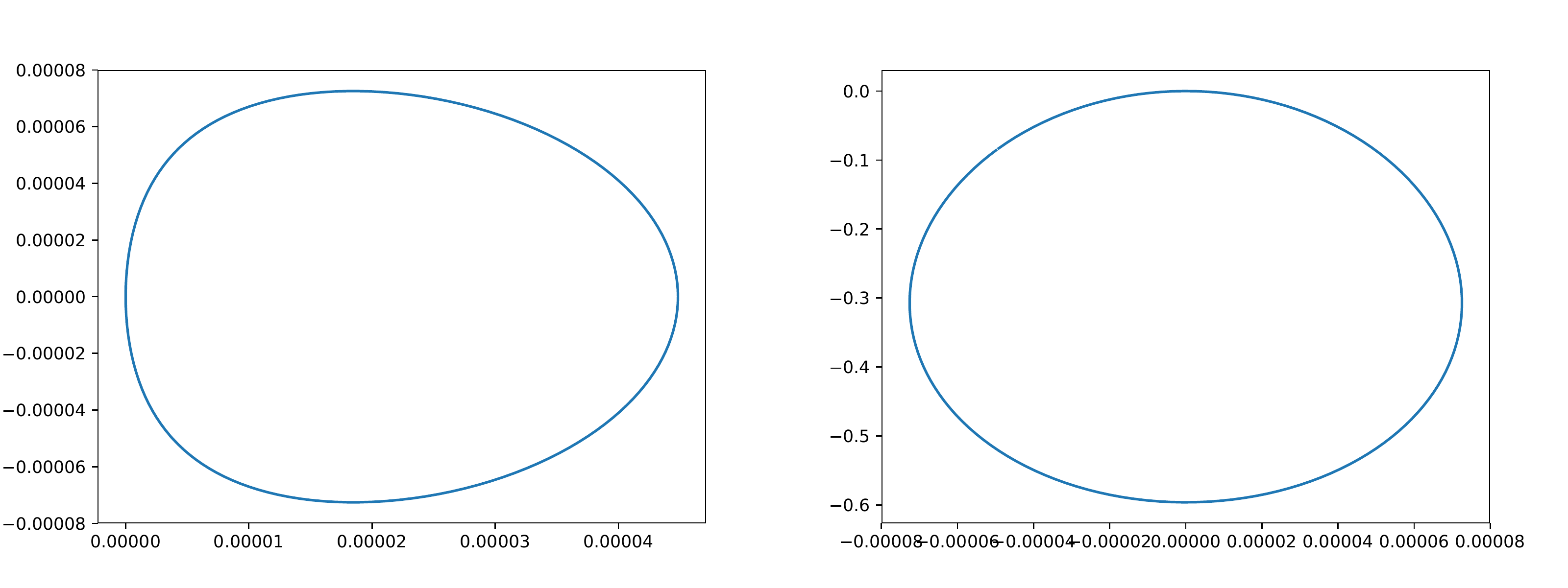}
  \caption{The $xy$- and the $yz$-projection of a periodic polar orbit in the lunar problem: $\mu=10^{-10}$ and $h=-1.5$ }
  \label{fig:lunar_c1_5}
\end{figure}

\begin{figure}[!htp]
  \centering
  \includegraphics[width=0.85\textwidth,clip]{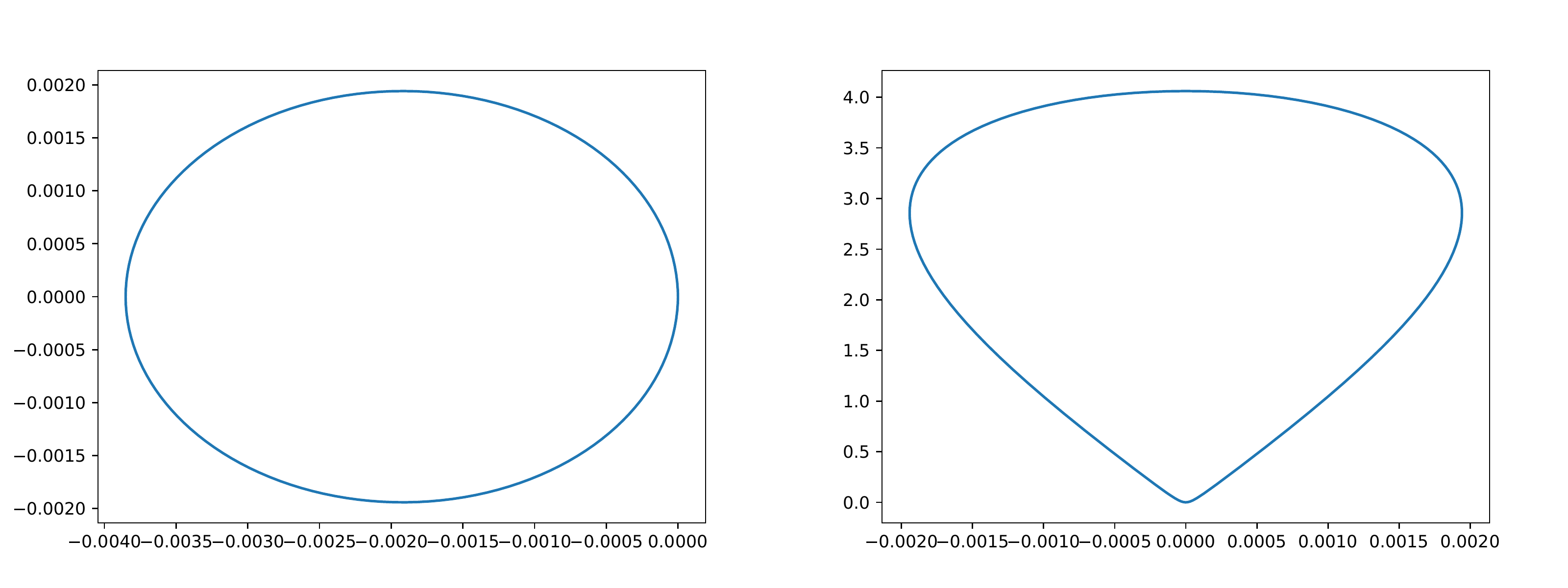}
  \caption{The $xy$- and the $yz$-projection of a periodic polar orbit in the lunar problem: $\mu=10^{-10}$ and $h=+8.0$ }
  \label{fig:lunar_c-8_0}
\end{figure}

Here we fix a small mass parameter $\mu>0$, namely $\mu=10^{-10}$ in the figures, and vary the energy. We plot the projection to the $xy$-plane and to the $yz$-plane, and look at the evolution as function of the energy. 
The typical situation for $h<0$ is drawn in Figure~\ref{fig:lunar_c1_5}.
For $h\geq 0$, the typical situation is drawn in \ref{fig:lunar_c-8_0} with the pointed tip in the $yz$-plane becoming more pronounced as the energy increases.

\begin{remark}
Periodic polar solutions can be found in Hill's lunar problem for all energies, even for $h>0$. Indeed, the Hill's region becomes unbounded, but remains bounded in the $z$-direction. This is of course not true for the restricted three-body problem, which we will discuss next.
\end{remark}

\subsubsection{Solutions for the restricted three-body problem}
The solutions we find for the rescaled Hamiltonian $H^\mu$ can be continued to larger $\mu$ as solutions for the unrescaled problem.
The energy is rescaled, and periodic polar orbits do not exist for all energies anymore. Indeed, for $h>0$ the orbits will typically escape a given region around the masses.

\begin{figure}[!htp]
  \centering
  \includegraphics[width=0.6\textwidth,clip]{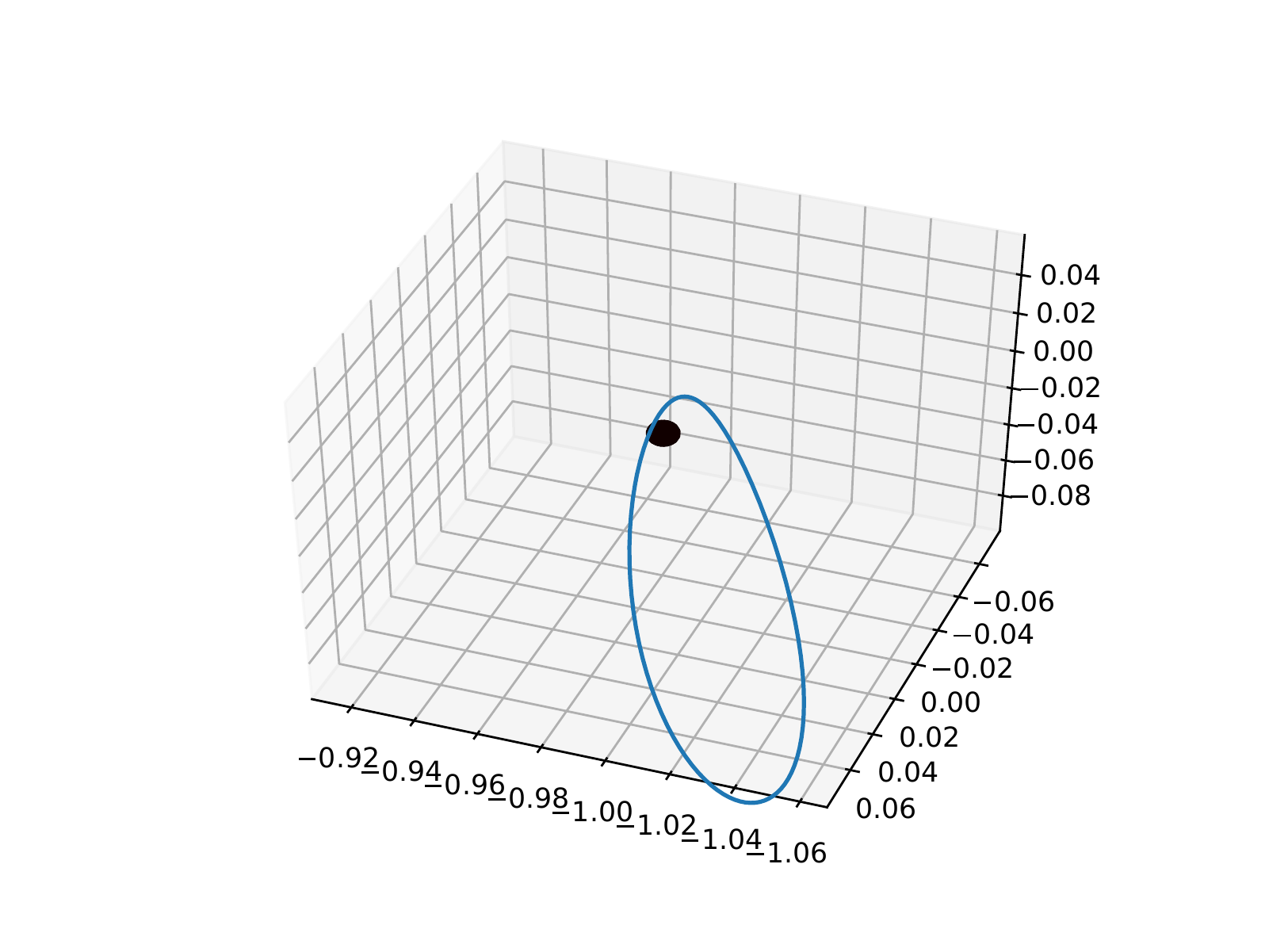}
  \caption{A periodic polar orbit in RTBP for the mass ratio Moon-Earth. The minimal distance to the center of the Moon, the black ball, is 4389 km. This orbit is mixed elliptic/negative hyperbolic.}
  \label{fig:moon_orbit}
\end{figure}

For small $\mu$ and suitably rescaled energy, the behavior of the orbits is of course the same as before. We will just mention one case that is of particular interest, namely the case $\mu\sim 0.01215$, which is the mass ratio of the Moon/Earth. We found that the periodic polar orbit can be continued to sufficiently large $h\sim -1.52$ such that it is no longer a physical collision orbit. This value of the Jacobi energy exceeds that of the first critical value. A 3d-plot of this orbit is given in Figure~\ref{fig:moon_orbit}.

We also remark that solutions in the restricted three-body problem include the families discussed in B81.
In the next section we see that the new family from the main theorem is for some energies a continuation of an orbit found in B81.
 
\subsubsection{Following the orbit for fixed energy and varying $\mu$}
Here we fix the Jacobi energy of the Hamiltonian $H_m$ and start at small, but positive $\mu$ at a near collision orbit.
We then homotope the mass parameter $\mu$ from this small value to $\mu=1$.
This gives part of a bridge connecting the polar orbits from this paper to the polar orbits from B81.
\begin{remark}
We point out that the part of the bridge indicated in the figures here is in the unrescaled problem. In other words, we are using the Hamiltonian $H_m$ rather than $H^\mu$.
\end{remark}
The bridge is constructed using a numerical homotopy. The homotopy becomes more difficult to carry out, i.e.~smaller parameter steps are needed, for larger $h$. In particular, for $h$ close to $-1.5$, the period of polar orbits near the smaller primary becomes very large.

\begin{figure}[!htp]
  \centering
  \includegraphics[width=0.85\textwidth,clip]{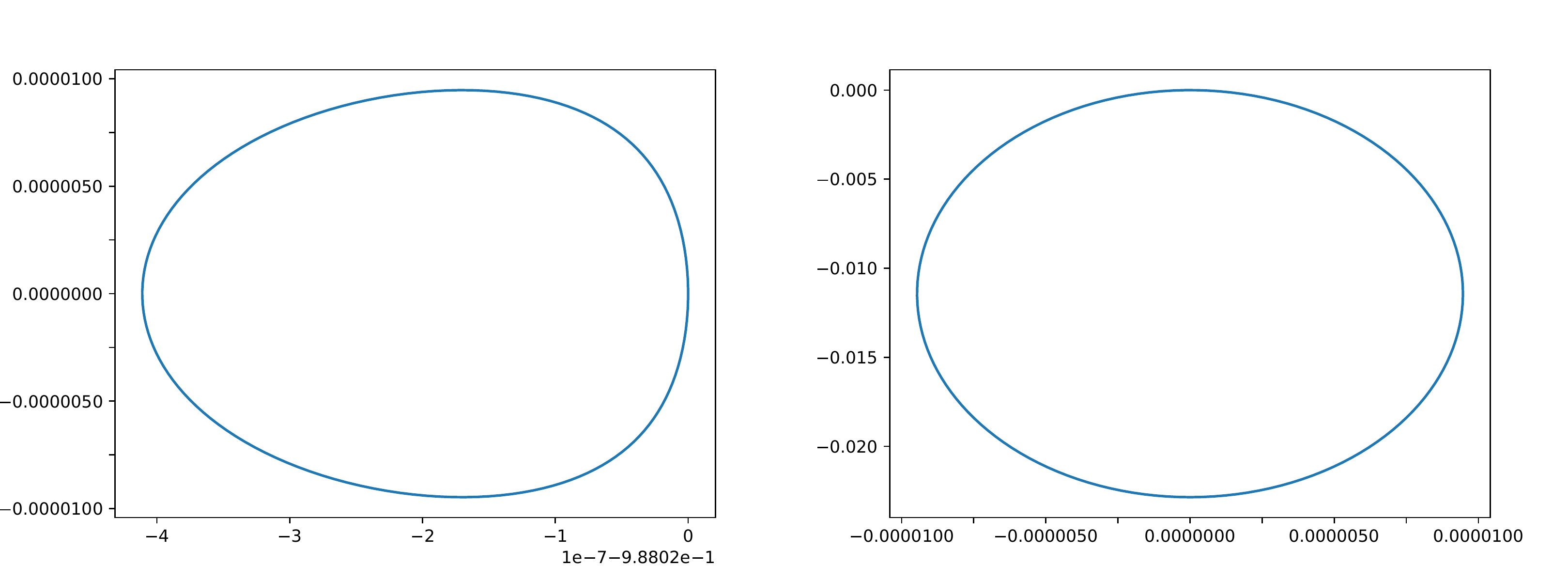}
  \caption{The $xy$- and $yz$-projection of an orbit in the bridge from lunar to rotating Kepler: $\mu=0.012$, $h=-2.0$.}
  \label{fig:bridge1}
\end{figure}

\begin{figure}[!htp]
  \centering
  \includegraphics[width=0.85\textwidth,clip]{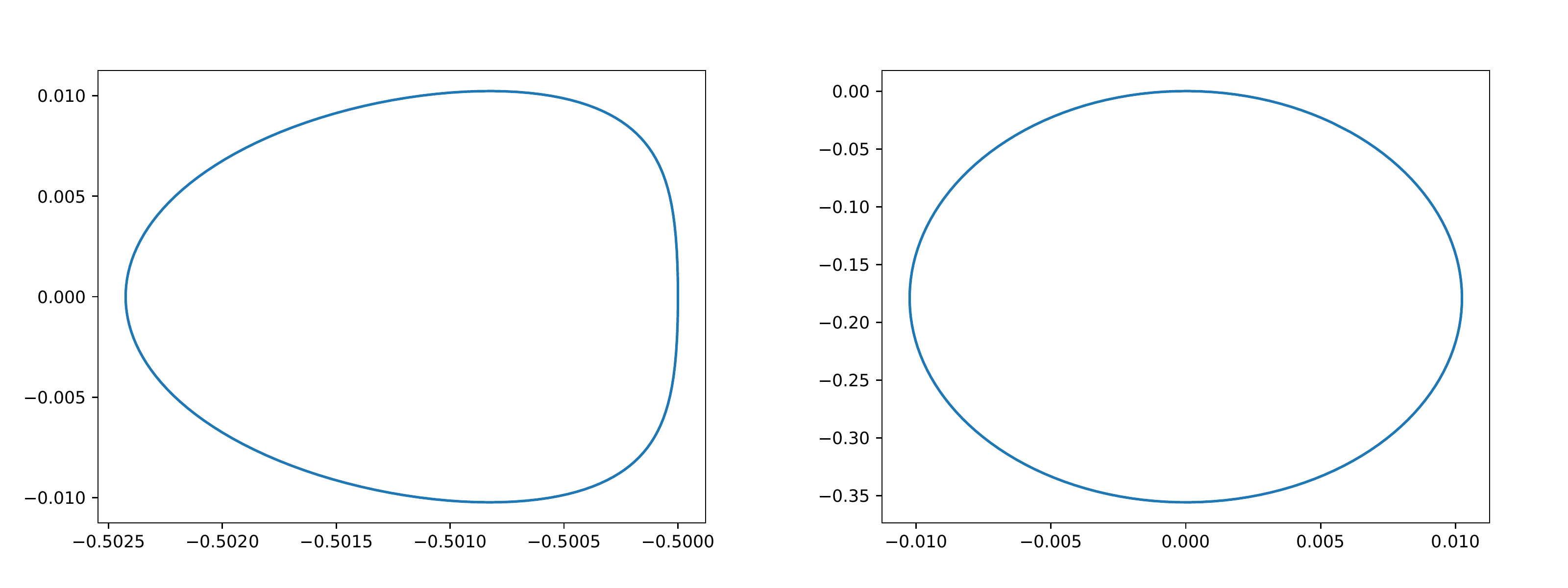}
  \caption{The $xy$- and $yz$-projection of an orbit in the bridge from lunar to rotating Kepler: $\mu=0.5$, $h=-2.0$.}
  \label{fig:bridge2}
\end{figure}

\begin{figure}[!htp]
  \centering
  \includegraphics[width=0.85\textwidth,clip]{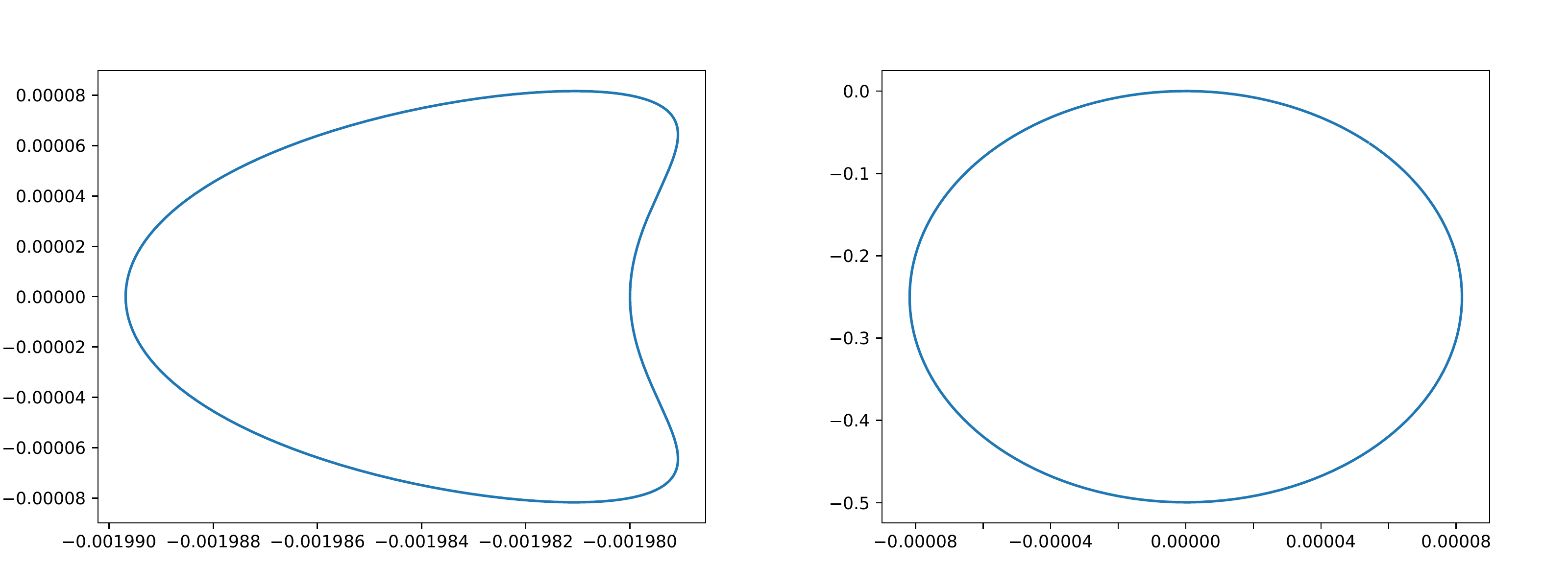}
  \caption{The $xy$- and $yz$-projection of an orbit in the bridge from lunar to rotating Kepler: $\mu=0.998$, $h=-2.0$.}
  \label{fig:bridge3}
\end{figure}

\subsubsection{Periapsis, apoapsis for Moon-Earth system}
We approximate the Moon-Earth system with the restricted three-body problem. For the distance Earth-Moon we take 386,000 km. Following the above scheme we find the polar orbit as a function of the Jacobi energy. It turns out that the polar orbit does not collide with the Moon for sufficiently large energy.
We have included a plot of the periapsis and of the apoapsis of the polar orbit.

\begin{figure}[!htp]
  \centering
  \includegraphics[width=0.5\textwidth,clip]{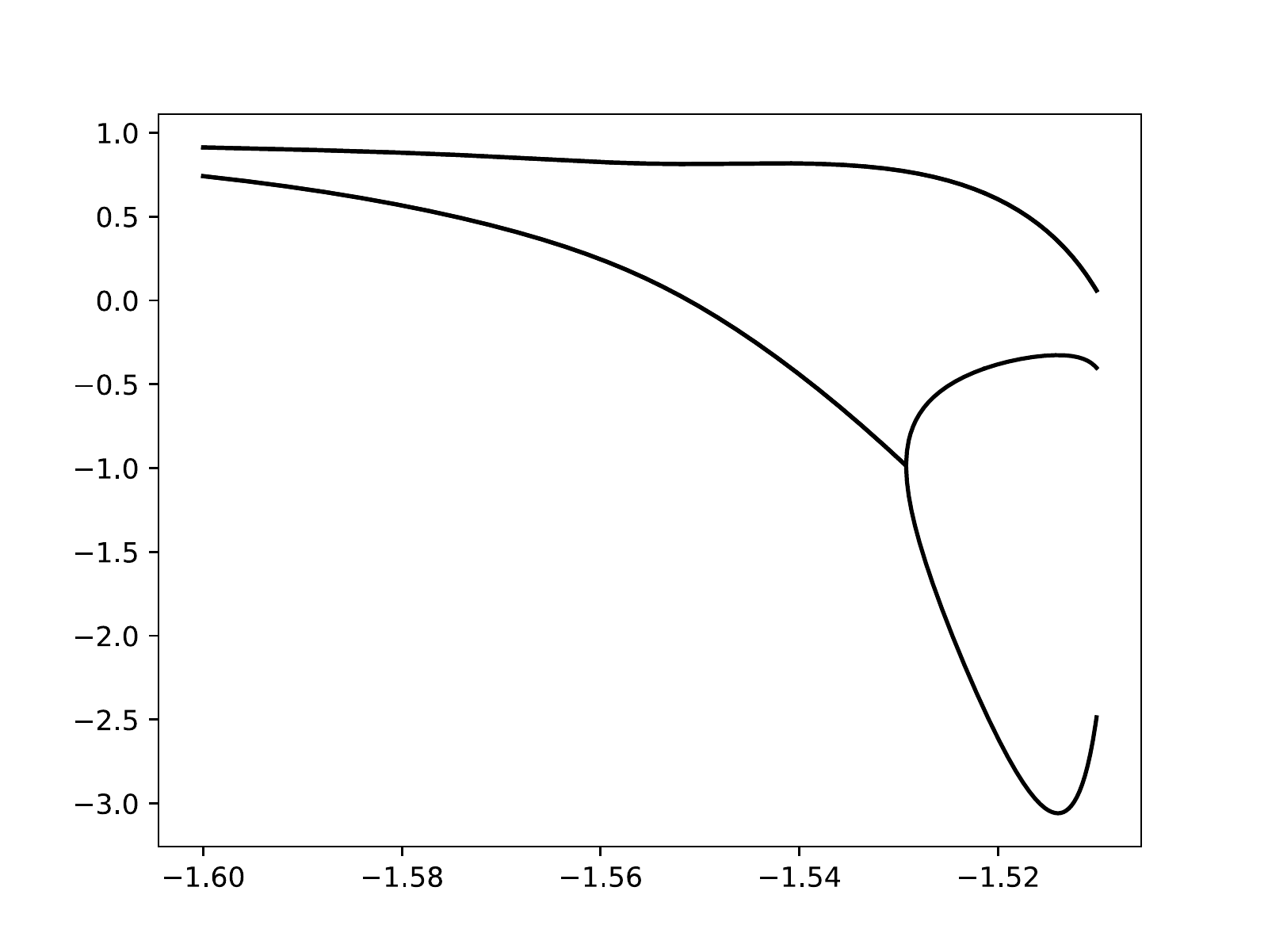}
  \caption{The real part of the eigenvalues of the linearized return map for the polar orbit in the Moon-Earth system as function of the Jacobi energy.}
  \label{fig:eigenvalues_moon}
\end{figure}

The stability properties of the polar orbit in the Moon-Earth system, though similar to those of the Hill's lunar system, are plotted in Figure~\ref{fig:eigenvalues_moon}. The bifurcation point indicates a period doubling/halving bifurcation.

\begin{figure}[!htp]
  \centering
  \includegraphics[width=1.0\textwidth,clip]{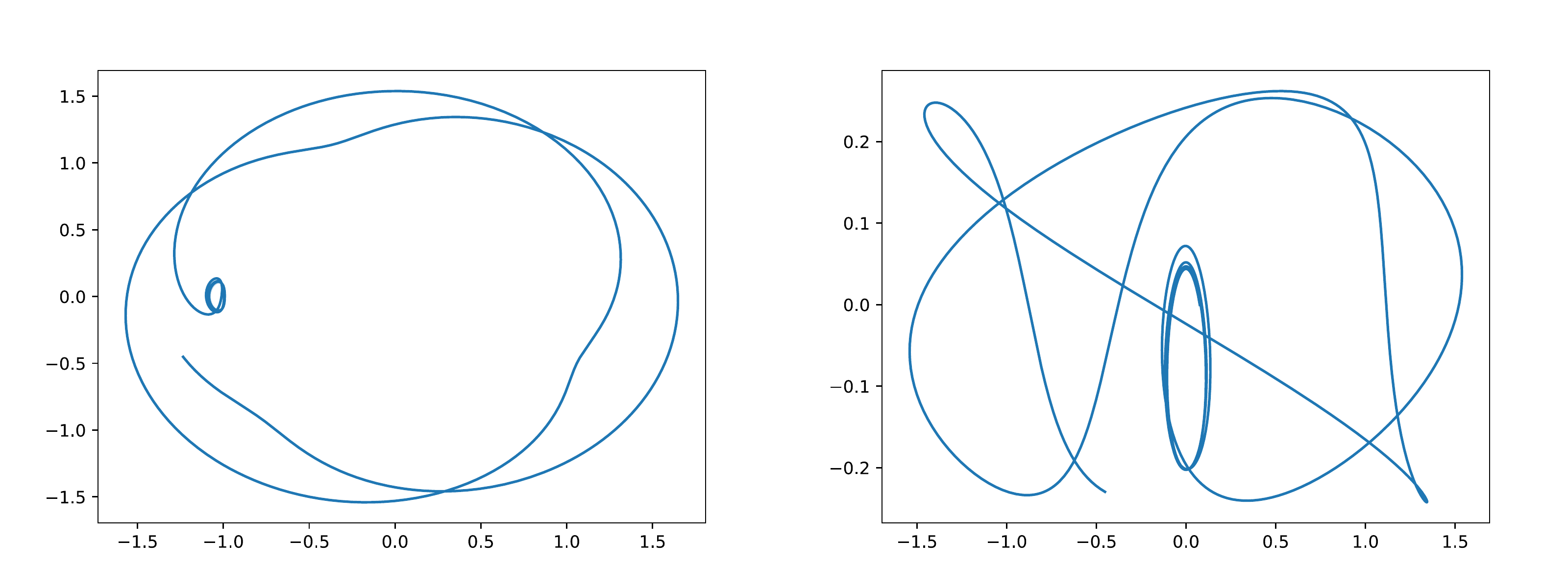}
  \caption{An orbit starting close to the periodic polar orbit: it shifted by about 400 km in the $y$-direction (orbits are fairly stable under shifts in the $x$-direction), and followed for five periods of the polar orbit.}
  \label{fig:unstable_moon}
\end{figure}

The effect of the instability that appears after the period doubling/halving bifurcation is indicated in Figure~\ref{fig:unstable_moon}.
We make the following observation: the periapsis of the polar orbit exceeds 1766 km (just 50km above the surface of the Moon) just before losing stability.
The values are so close though that the approximations we made (circular restricted three-body problem) most likely spoils stability of a physical non-collision orbit.

\section{Geometric and Global Properties}
\label{sec:Geometry} 
In this section we will give a theoretical overview of results that explain the bridge. These results are limited to bounded Hill's region, meaning that the Jacobi parameter should be chosen to be sufficiently negative; for example $h=-2.0$ will certainly work. The numerical results go beyond this, namely the bridge exists for much higher Jacobi energy $h$, but the  theoretical underpinning here does not directly apply.

Recall that $-\tfrac{3^{4/3}}{2}$ is the critical value of Hill's lunar problem. Pick $h_0<-\tfrac{3^{4/3}}{2}$ an energy
below the critical value for which the consecutive collision orbit on the $z$-axis in Hill's lunar problem is non-degenerate.
Choose
$$h \colon [0,1] \to \mathbb{R}$$
a smooth function with the property that 
$$h(\mu)<h_1(\mu)$$
for every $\mu \in [0,1]$, where $h_1(\mu)$ is the first critical value of the restricted three-body problem with mass
parameter $\mu$, and moreover for small $\mu$, $h$ is given by
$$h(\mu)=-\tfrac{3}{2}+c_0 \mu^{2/3}.$$
For $\mu_\infty \in (0,1]$ and $\mu \in [0,\mu_\infty)$ let $v_\mu$ be the smooth family of symmetric periodic orbits 
such that $v_0$ is the consecutive collision orbit in Hill's lunar problem of energy $h_0$ and $v_\mu$
is a symmetric periodic orbit in the restricted three-body problem of energy $h(\mu)$ and mass parameter $\mu$ for every
$\mu \in (0,\mu_\infty)$. Because $v_0$ is non-degenerate such a family exists uniquely by the implicit function theorem
and moreover we can assume that $v_\mu$ is non-degenerate for every $\mu \in [0,\mu_\infty)$ and $\mu_\infty$
is maximal with this property. Here we consider the regularized restricted three-body problem and it might happen that
$v_\mu$ has collisions. We denote by $\Omega$ the $\omega$-limit set of the family $v_\mu$ consisting of all
symmetric periodic orbits $w$ for which there exists a sequence $\mu_\nu \in [0,\mu_\infty)$ converging to $\mu_\infty$ 
such that $w=\lim_{\nu \to \infty} v_{\mu_\nu}$. In \cite{albers-frauenfelder-koert-paternain} it was shown that
the bounded components of the regularized energy hypersurfaces in the planar restricted three-body problem are of contact type for energies below and slightly above
the first critical value. Cho and Kim have generalized this result to the spatial case, \cite{ChoKim}. Using this result and the results in \cite{belbruno-frauenfelder-vankoert} we obtain the following information about $\Omega$
\begin{theorem}
If $\mu_\infty<1$, then $\Omega$ is nonempty, compact and connected and consists of degenerate symmetric periodic orbits
in the restricted three-body problem of mass parameter $\mu_\infty$ and energy $h(\mu_\infty)$. 
\end{theorem}
{\em Proof sketch:} The proof of the Theorem relies on the Theorem of Arzela-Ascoli. Because $h(\mu)<h_1(\mu)$ the Hill's region
is bounded and therefore the images of $v_\mu$ lie in a compact set. The contact condition guarantees that the periods of
the periodic orbits are uniformly bounded from which equicontinuity follows. Details can be found in \cite{belbruno-frauenfelder-vankoert}.

Using local Rabinowitz Floer homology we obtain the following result \cite{belbruno-frauenfelder-vankoert}.
\begin{theorem}
Assume that $\Omega$ is isolated and for an open neighborhood $\mathcal{U}$ of $\Omega$ there exists a sequence
$\mu_\nu \in (\mu_\infty,1]$ converging to $\mu_\infty$ with the property that there are no symmetric periodic orbits
in $\mathcal{U}$ for mass parameter $\mu_\nu$ and energy $h(\mu_\nu)$. Then there exists $\mu_0 \in [0,\mu_\infty)$ such that
for every $\mu \in [\mu_0,\mu_\infty)$ there exists in $\mathcal{U}$ a symmetric periodic orbit for mass parameter $\mu$
and energy $h(\mu)$ different from $v_\mu$.
\end{theorem}

\begin{remark}
We briefly remark that without the contact condition families of periodic orbits can end in a blue sky catastrophe, meaning that the period goes to infinity.
\end{remark}

\ack

Edward Belbruno would like to acknowledge the support of Humboldt Stiftung of the Federal Republic of Germany that made this research possible and the support of the University of Augsburg for his visit from 2018-19. Urs Frauenfelder was supported by DFG grant FR 2637/2-1, of the German government.
Otto van Koert was supported by NRF grant NRF-2016R1C1B2007662, funded by the Korean Government.

\appendix

\subsection{Regularization in coordinates}
Moser regularization is based on stereographic projection. The basic idea is to switch the role of momentum and position in the $q,p$-coordinates, and use the $p$-coordinates as position coordinates in $T^*\mathbb{R}^n$, where we think of $\mathbb{R}^n$ as a chart for $S^n$.
To prevent confusion, we set $\vec x=p$ and $\vec y=q$.

To go from $T^*S^n$ to $T^*\mathbb{R}^n$ we use the map
\begin{equation}
\label{eq:regularized_to_unregularized}
\begin{split}
\vec x &= \frac{\vec \xi}{1-\xi_0} \\
\vec y &= \eta_0 \vec \xi +(1-\xi_0) \vec \eta
\end{split}
\end{equation}
We think of $T^*S^n$ as the following submanifold of $T^*\mathbb{R}^{n+1}$.

To go from $T^*\mathbb{R}^n$ to $T^*S^n$, we use the inverse given by
\begin{equation}
\label{eq:unregularized_to_regularized}
\begin{split}
\xi_0 &= \frac{|\vec x|^2-1}{|\vec x|^2+1} \\
\vec \xi &= -\frac{2 \vec x}{|\vec x|^2+1} \\
\eta_0 &= -\langle \vec x,\vec y \rangle \\
\vec \eta &= \frac{|\vec x|^2+1}{2}\vec y 
- \langle \vec x, \vec y \rangle \vec x
\end{split}
\end{equation}

The Belbruno transform employs a M\"obius transformation which sends to the collision point $\infty$ to $1=(1,0,\ldots,0)\in \mathbb{R}^n$.
In coordinates (the index $j=2,3$), the forward Belbruno transformation is given by
\[
\begin{split}
Q_1 &= \frac{1-\Vert p\Vert^2}{2}q_1 +\langle q, p\rangle (p_1+1) \\
Q_j & = \frac{\Vert p\Vert^2+1}{2}q_j + p_1 q_j - p_j q_1 - \langle q, p\rangle p_j\\
P_1 &= \frac{\Vert p\Vert^2 -1 }{\Vert p+1 \Vert^2} \\
P_j &= \frac{2p_j}{\Vert p+1 \Vert^2}
\end{split}
\]
The inverse Belbruno transform is given by
\[
\begin{split}
    q_1 &= \frac{1-\Vert P\Vert^2}{2}Q_1 + \langle Q,P\rangle (P_1 - 1)\\
    q_j &= \frac{\Vert P\Vert^2}{2}Q_j - P_1 Q_j + P_j Q_1 - \langle Q,P\rangle P_j\\
    p_1 &= \frac{1-\Vert P\Vert^2}{ \Vert P-1 \Vert^2 } \\
    p_j &= \frac{2P_j}{ \Vert P-1 \Vert^2 } 
\end{split}
\]

\subsection{Hamiltonian vector field with constraints}
The setup is the following. We are given a manifold $M$, which is a symplectic submanifold of the symplectic manifold $(N,\Omega)$. We denote the inclusion by $\iota: M\to N$, and the induced symplectic form on $M$ by $\omega:=\iota^*\Omega$. We assume that $M=f_1^{-1}(0) \cap f_2^{-1}(0)$.
In addition, we are given a Hamiltonian function $H_N:N \to \mathbb{R}$, and we have the induced Hamiltonian $H_M=\iota^*H_N$.
In our case $N=T^*\mathbb{R}^{n+1}$ and 
$$
M:=T^*S^n
=\{ (\xi,\eta) \in T^*\mathbb{R}^{n+1} ~|~\Vert\xi \Vert^2=1,\quad 
\langle \xi,\eta\rangle =0
\}
.
$$ 
The functions that define $M$ are
$$
f_1=\frac{1}{2}\Vert\xi \Vert^2-\frac{1}{2},\quad
f_2=\langle \xi,\eta\rangle
.
$$
In our case, the symplectic manifold $N=T^*\mathbb{R}^{n+1}$ has a global chart, but $T^*S^n$  has not. We will give a formula for the Hamiltonian vector field $X_H$ on $M$ in terms of Hamiltonian vector field on $N$. In our example, this means that we can use the global coordinates on $N=T^*\mathbb{R}^{n+1}$.
We have
\begin{equation}
\label{eq:contrained_ham_vf}
X_H=X_{H_N}+c_1 X_{f_1}+c_2 X_{f_2}
\end{equation}
where
\[
\begin{split}
X_{H_N}&=\sum_{j=0}^n \frac{\partial H_N}{\partial \eta_j} \frac{\partial}{\partial \xi_j}-\frac{\partial H_N}{\partial \xi_j} \frac{\partial}{\partial \eta_j} \\
X_{f_1}&=\sum_{j=0}^n \frac{\partial f_1}{\partial \eta_j} \frac{\partial}{\partial \xi_j}-\frac{\partial f_1}{\partial \xi_j} \frac{\partial}{\partial \eta_j} = -\sum_j \xi_j \frac{\partial}{\partial \eta_j}\\
X_{f_2}&=\sum_{j=0}^n \frac{\partial f_2}{\partial \eta_j} \frac{\partial}{\partial \xi_j}-\frac{\partial f_2}{\partial \xi_j} \frac{\partial}{\partial \eta_j} = \sum_j \xi_j \frac{\partial}{\partial \xi_j}-\eta_j \frac{\partial}{\partial \eta_j}\\
c_1 &= \frac{df_2(X_{H_N})}{df_2(X_{f_1})} = -\frac{ \{ f_2,H_N\} }{\{ f_1,f_2\}}
=
-\{ f_2, H_N\}
\\
c_2 &= \frac{df_1(X_{H_N})}{df_1(X_{f_2})} = \frac{ \{ f_1,H_N\} }{\{ f_1,f_2\}}
=
\{ f_1, H_N\}
\end{split}
\]
The Poisson brackets, defined by $\{ f, g\}:=\omega(X_f, X_g\}$ are of course not needed to do the computations, but the clarify that the situation is if $M$ is symplectic submanifold of higher codimension, where a matrix filled with $\{ f_i, f_j \}$ has to be inverted.
A computation shows that the above vector field is tangent to the submanifold $M$ and that it is the Hamiltonian vector field.

\end{document}